\documentstyle[12pt]{article}
\input amssym.def
\begin{document}
\def \Z{\Bbb Z}
\def \C{\Bbb C}
\def \R{\Bbb R}
\def \Q{\Bbb Q}
\def \N{\Bbb N}
\def \wt{{\rm wt}}
\def \tr{{\rm tr}}
\def \span{{\rm span}}
\def \Res{{\rm Res}}
\def \End{{\rm End}}
\def \Ind {{\rm Ind}}
\def \Irr {{\rm Irr}}
\def \Aut{{\rm Aut}}
\def \Hom{{\rm Hom}}
\def \mod{{\rm mod}}
\def \ann{{\rm Ann}}
\def \<{\langle} 
\def \>{\rangle} 
\def \t{\tau }
\def \a{\alpha }
\def \e{\epsilon }
\def \l{\lambda }
\def \L{\Lambda }
\def \g{\gamma}
\def \b{\beta }
\def \om{\omega }
\def \o{\omega }
\def \c{\chi}
\def \ch{\chi}
\def \cg{\chi_g}
\def \ag{\alpha_g}
\def \ah{\alpha_h}
\def \ph{\psi_h}

\def \bconj{\begin{conj}\label}
\def \econj{\end{conj}}
\def \be{\begin{equation}\label}
\def \ee{\end{equation}}
\def \bex{\begin{exa}\label}
\def \eex{\end{exa}}
\def \bl{\begin{lem}\label}
\def \el{\end{lem}}
\def \bt{\begin{thm}\label}
\def \et{\end{thm}}
\def \bp{\begin{prop}\label}
\def \ep{\end{prop}}
\def \br{\begin{rem}\label}
\def \er{\end{rem}}
\def \bc{\begin{coro}\label}
\def \ec{\end{coro}}
\def \bd{\begin{de}\label}
\def \ed{\end{de}}
\def \pf{{\bf Proof. }}
\def \voa{{vertex operator algebra}}

\newtheorem{thm}{Theorem}[section]
\newtheorem{prop}[thm]{Proposition}
\newtheorem{coro}[thm]{Corollary}
\newtheorem{conj}[thm]{Conjecture}
\newtheorem{exa}[thm]{Example}
\newtheorem{lem}[thm]{Lemma}
\newtheorem{rem}[thm]{Remark}
\newtheorem{de}[thm]{Definition}
\newtheorem{hy}[thm]{Hypothesis}
\makeatletter
\@addtoreset{equation}{section}
\def\theequation{\thesection.\arabic{equation}}
\makeatother
\makeatletter

\newcommand{\rw}{\rightarrow}
\newcommand{\n}{\:^{\circ}_{\circ}\:}

\begin{center}{\Large \bf Regular representations of 
vertex operator algebras, I}
\end{center}

\begin{center}{Haisheng Li\footnote{Supported by NSF grant 
DMS-9616630}\\
Department of Mathematical Sciences\\
Rutgers University-Camden\\
Camden, NJ 08102}
\end{center}
\section{Introduction}

In a previous study [Li3], the physical super selection
principle in vertex operator algebra theory
seems to tell us that for a vertex operator algebra $V$,
one should be able to
obtain each irreducible $V$-module from the adjoint module $V$ by
changing certain things. To a certain extent, what we are expecting 
is an analogue of the
orbital construction of representations for a Lie group [Ki].
This was our original motivation for this study on
regular representations of vertex operator algebras.

This is the first paper in a series for this study.
In this paper, 
given a module $W$ for a vertex operator algebra $V$ 
and a nonzero complex number $z$ we construct 
a canonical (weak) $V\otimes V$-module ${\cal{D}}_{P(z)}(W)$
(a subspace of $W^{*}$ depending on $z$). We prove that for $V$-modules 
$W, W_{1}$ and $W_{2}$, a $P(z)$-intertwining map
of type ${W'\choose W_{1}W_{2}}$ ([H3], [HL0-3]) exactly 
amounts to a $V\otimes V$-homomorphism from 
$W_{1}\otimes W_{2}$ into ${\cal{D}}_{P(z)}(W)$.
Using Huang and Lepowsky's one-to-one linear correspondence between 
the space of intertwining operators and the space of $P(z)$-intertwining 
maps of the same type we obtain a canonical linear isomorphism from
the space ${\cal{V}}^{W'}_{W_{1}W_{2}}$ of
intertwining operators of the indicated type to
$\Hom _{V\otimes V}(W_{1}\otimes W_{2},{\cal{D}}_{P(z)}(W))$.
In the case that $W=V$, we obtain a decomposition of
Peter-Weyl type for ${\cal{D}}_{P(z)}(V)$,
which are what we call the regular representations of $V$. 

Let us start with the classical Peter-Weyl theory.
Let $G$ be a group $G$ and $F(G,{\C})$ be the space of
complex-valued functions on $G$. Then $F(G,{\C})$
is a $G\times G$-module with
\begin{eqnarray}\label{egbim}
((g_{1},g_{2})f)(g)=f(g_{1}^{-1}gg_{2})
\end{eqnarray}
for $g_{1}, g_{2}, g\in G, f\in F(G,{\C})$.
Furthermore, if $G$ is a topological (Lie, algebraic) group,
all continuous (${\C}^{\infty}$, algebraic) functions
form a $G\times G$-submodule.
Now let $G$ be a compact Lie group and $C^{0}(G)$ be the space of 
continuous functions on $G$. A representative function on $G$ is
a continuous function which generates a finite-dimensional
$1\otimes G$-submodule of $C^{0}(G)$, and 
the regular representation of $G$ is the $G\times G$-module
${\cal{J}}(G,{\C})$ of representative functions on $G$ (cf. [BD]).  
Let $U$ be an irreducible $G$-module. For $u\in U,\; u^{*}\in U^{*}$, 
we can view $u^{*}\otimes u$ as a function on $G$ by
$$(u^{*}\otimes u)(g)=\<u^{*},gu\>\;\;\;\mbox{ for }g\in G.$$
Then $U^{*}\otimes U$ is embedded into $C^{0}(G)$ as a $G\otimes G$-module.
A theorem of Peter-Weyl type states that
${\cal{J}}(G,{\C})$ is a dense subspace of $C^{0}(G)$
and that ${\cal{J}}(G,{\C})$ is a 
direct sum of $U^{*}\otimes U$, where $U$ runs through a complete  set
of representatives of equivalence classes of finite-dimensional 
irreducible $G$-modules. 

To start vertex operator algebra theory, let us recall from
[FLM] and [FHL] (see also [B]) the
definition of vertex operator algebra, which will be the official
definition for this paper.
{\em A vertex operator algebra} is a ${\Z}$-graded vector
space $V=\coprod_{n\in {\Z}}V_{(n)}$ such that $\dim V_{(n)}<\infty$ 
for all $n\in {\Z}$ and $V_{(n)}=0$ for $n$ sufficiently small, 
equipped with a bilinear ``vertex multiplication''
\begin{eqnarray}
Y(\cdot,x)\cdot: & &V\otimes V\rightarrow V((x))\nonumber\\ 
& &u\otimes v\mapsto Y(u,x)v=\sum_{n\in {\Z}}u_{n}vx^{-n-1}
\end{eqnarray}
such that for $m, n, k\in {\Z},\;u\in V_{(k)}$,
\begin{eqnarray}\label{egrad}
u_{n}V_{(m)}\subset V_{(m+k-n-1)},
\end{eqnarray}
and such that the {\em Jacobi identity} holds for $u,v\in V$:
\begin{eqnarray}\label{e0jacobi}
& &x_{0}^{-1}\delta\left(\frac{x_{1}-x_{2}}{x_{0}}\right)
Y(u,x_{1})Y(v,x_{2})
-x_{0}^{-1}\delta\left(\frac{x_{2}-x_{1}}{-x_{0}}\right)
Y(v,x_{2})Y(u,x_{1})\nonumber\\
& &=x_{2}^{-1}\delta\left(\frac{x_{1}-x_{0}}{x_{2}}\right)
Y(Y(u,x_{0})v,x_{2}).
\end{eqnarray}
It is equipped with a vector ${\bf 1}$, 
called the {\em vacuum vector},
such that for $v\in V$,
\begin{eqnarray}
& &Y({\bf 1},x)v=v,\\
& &Y(v,x){\bf 1}\in V[[x]]\;\;\;\mbox{ and  }
\;\lim_{x\rightarrow 0}Y(v,x){\bf 1}=v.
\end{eqnarray}
It is also equipped with a vector $\omega$, called the 
{\em Virasoro element}, such that
\begin{eqnarray}
[L(m), L(n)]=(m-n)L(m+n)+\frac{m^{3}-m}{12}\delta_{m+n,0}(\mbox{rank}V)
\end{eqnarray}
for $m,n\in {\Z}$, where $Y(\omega,x)=\sum_{n\in {\Z}}L(n)x^{-n-2}$ and
$\mbox{rank}V$ is a complex number, called the {\em rank} of $V$, and 
such that for $v\in V_{(n)}, n\in {\Z}$,
\begin{eqnarray}
& &Y(L(-1)v,x)={d\over dx}Y(v,x),\\
& & L(0)v=nv.
\end{eqnarray}

This completes the definition. It is well known ([FLM], [FHL]) that 
the Jacobi identity is 
equivalent to the following rationality, 
commutativity and associativity:

(rationality) For $u,v,w\in V,\; w'\in V'$, where 
$V'=\coprod_{n}V_{(n)}^{*}$, the formal series
\begin{eqnarray}\label{e02}
\<w', Y(u,x_{1})Y(v,x_{2})w\>\left(=\sum_{m,n\in {\Z}}\<w',u_{m}v_{n}w\>
x_{1}^{m-1}x_{2}^{-n-1}\right)
\end{eqnarray}
absolutely converges to a rational function of the form
\begin{eqnarray}\label{ef}
f(x_{1},x_{2})=\frac{g(x_{1},x_{2})}{(x_{1}-x_{2})^{k}}
\end{eqnarray}
in the domain $|x_{1}|>|x_{2}|>0$, where 
$g\in {\C}[x_{1},x_{1}^{-1},x_{2},x_{2}^{-1}]$ and $k$ is an integer
only depending on $u$ and $v$; independent of $w'$ and $w$; 

(commutativity) The formal series $\<w', Y(v,x_{2})Y(u,x_{1})w\>$
absolutely converges to the same rational function $f(x_{1},x_{2})$
in the domain $|x_{2}|>|x_{1}|>0$;
 
(associativity) The formal series
\begin{eqnarray}\label{e03}
\<w', Y(Y(u,x_{0})v,x_{2})w\>\left(=\sum_{m,n\in {\Z}}\<w', (u_{m}v)_{n}w\>
x_{0}^{-m-1}x_{2}^{-n-1}\right)
\end{eqnarray}
absolutely converges to the rational function $f(x_{0}+x_{2},x_{2})$
in the domain $|x_{2}|>|x_{0}|>0$.
If we specialize $x_{2}$ to a nonzero complex number $z$, (\ref{e02})
gives rise to a meromorphic 
function on the sphere $CP^{1}={\C}\cup \{\infty\}$
with only three possible poles at $0, \infty$ and $z$ while
(\ref{e03}) gives rise to a meromorphic 
function on the sphere $CP^{1}$
with only three possible poles at $0, \infty$ and $-z$. 
The rationality, commutativity and associativity
is the basis for Huang's geometric interpretation of a vertex 
operator algebra [H1].

For a vertex operator algebra $V$, we have the following fundamental 
results established in [FHL]. It was proved that 
$V\otimes V$ has a natural vertex operator algebra structure and that
for any $V$-modules $W_{1}$ and $W_{2}$, $W_{1}\otimes W_{2}$ is a
natural $V\otimes V$-module; 
Let $W=\coprod_{h\in {\C}}W_{(h)}$ be a $V$-module and
let $W'=\coprod_{h\in {\C}}W_{(h)}^{*}$, the restricted dual.
For $v\in V$, define $Y'(v,x)$ to be an element of
$({\rm End}\;W')[[x,x^{-1}]]$ by 
\begin{eqnarray}
\<Y'(v,x)w',w\>=\<w',Y(e^{xL(1)}(-x^{-2})^{L(0)}v,x^{-1})w\>
\end{eqnarray}
for $w'\in W',\; w\in W$.
Then $(W', Y')$ carries the structure of a $V$-module; For $v\in V$, set
\begin{eqnarray}\label{eyo}
Y^{o}(v,x)=Y(e^{xL(1)}(-x^{-2})^{L(0)}v,x^{-1}).
\end{eqnarray}
It was proved ([HL1], [FHL])  that
$Y^{0}$ satisfies the opposite Jacobi identity (\ref{ehl}), so that the pair
$(W,Y^{o})$ carries the structure of a right $V$-module.
These facts should be compared with those
mentioned in the second paragraph for a group.

Now we consider what the space of a regular representation should be. 
The classical analogue suggests that we consider ``functions on $V$
of the form $w'\otimes w$'' for $w'\in W',\; w\in W$, where
$W$ is an irreducible $V$-module. By analogy
we should view $w'\otimes w$ as
a ${\C}[x,x^{-1}]$-valued linear function on $V$
by considering $\<w',Y(v,x)w\>$ for $v\in V$.
To get a ${\C}$-valued function on $V$,
we use Huang and Lepowsky's idea to 
evaluate the formal variable $x$ at a nonzero complex 
number $z$ as in [H3] and [HL0-3], so that $w'\otimes w$ 
can be viewed as a linear functional on $V$ by
\begin{eqnarray}
(w'\otimes w)(v)=\<w',Y(v,z)w\>\;\;\;\mbox{ for }v\in V.
\end{eqnarray}
Such functionals satisfy 
a certain nice rationality property as described below.
For a linear functional $\alpha$ on $V$, we naturally extend
$\alpha$ to a linear function on $V[[x,x^{-1}]]$.
Then consider the formal series $\alpha(Y(u,x)v)$ with $u,v\in V$.
If $\alpha=w'\otimes w$, we have
$$\alpha(Y(u,x)v)=\<w',Y(Y(u,x)v,z)w\>,$$
which by associativity converges to a rational function 
in the domain $0<|x|<|z|$ with only two possible (finite) poles at $0$ 
and $-z$. Globally, from $Y$ we obtain a linear map $F$ from 
$W'\otimes W$ to $V^{*}$ defined by
\begin{eqnarray}
\<F(w'\otimes w), v\>=\<w',Y(v,z)w\>
\end{eqnarray}
for $w'\in W',\; w\in V,\; v\in V$. This map $F$
is nothing but a $Q(z)$-intertwining map as known in [HL1].
Linear functionals $\alpha$ with this
rationality property are called $Q(z)$-linear functions and 
will be studied in [Li6].

In this paper we shall study the so-called $P(z)$-linear 
functionals defined below. Notice that previously we basically used 
the canonical intertwining operator $Y$ of 
type ${W\choose VW}$. It was known ([FHL], [HL2]) 
that there are canonical linear isomorphisms from the space of
intertwining operators of type ${W\choose VW}$ to
 the space of intertwining operators of type ${V'\choose W'W}$.
Let ${\cal{Y}}$ be an intertwining operator of type ${V'\choose W'W}$.
(One can make use of the vacuum vector to get a canonical one.)
Then ${\cal{Y}}(w',z)w$ is a linear functional on $V$.
(To be rigorous we have to choose a branch of log function 
to evaluate $x$ at $z$ for ${\cal{Y}}(w',x)w$, see Section 4, 
or [HL1]. Here we neglect this issue for this introduction.)
This gives another way to 
view $w'\otimes w$ as a linear functional on $V$.
This time, for $u,v\in V$, the formal series
$$\<Y(u,x){\cal{Y}}(w',z)w,v\>\;\;\;\;
\left(=\< {\cal{Y}}(w',z)w,Y^{o}(u,x)v\>\right)$$
by commutativity converges to a rational function in 
the domain $|x|>|z|$ with only two possible poles
at $0$ and $z$, or a meromorphic function on $CP^{1}$
with only three possible poles at $0, z$ and $\infty$.
Globally, evaluated intertwining operators ${\cal{Y}}(\cdot,z)\cdot$ 
at $x=z$ are exactly captured by the notion of $P(z)$-intertwining 
map ([H3], [HL0-3]).
Then we consider linear functionals $\alpha$ on $V$ such that
 for $u,v\in V$, the formal series
$\alpha(u,v,x):=\alpha(Y^{o}(u,x)v)$ converges to a rational 
function of $x$ with only two possible poles at $0$ and $z$, and
we call such linear functionals $P(z)$-{\em linear 
functionals} on $V$ and denote by ${\cal{D}}_{P(z)}(V)$ the space of 
all $P(z)$-linear functionals. 
We define $Y^{R}(u,x)\alpha$ and $Y^{L}(u,x)\alpha$ by requiring 
$\< Y^{R}(u,x)\alpha,v\>$ to be the formal Laurent series expansion 
of the rational function $\alpha(u,v,x)$ in the domain $0<|x|<|z|$, 
and $\< Y^{L}(u,x)\alpha,v\>$ to be the Laurent series expansion 
of the rational function $\alpha(u,v,x+z)$ in the domain $0<|x|<|z|$.
Our first main result is that
$Y^{L}$ and $Y^{R}$ give rise to
a structure of a weak $V\otimes V$-module on
${\cal{D}}_{P(z)}(V)$ in the sense that all the axioms defining
the notion of a module except those involving the grading hold.

Vertex operator algebras are similar to Lie algebras.
One of the similarities is that ideals are always two-sided
due to the skew-symmetry. Consequently,
a simple vertex operator algebra 
itself is an irreducible module. Most of the time, the adjoint module
and other modules can be equally treated.
Here, we construct a weak $V\otimes V$-module ${\cal{D}}_{P(z)}(W)$ 
for a general $V$-module $W$. 
Then, for any additional modules $W_{1}$ and $W_{2}$, we 
identify the space ${\cal{M}}[P(z)]^{W'}_{W_{1}W_{2}}$ of 
$P(z)$-intertwining maps of the indicated type with 
 ${\rm Hom}_{V\otimes V}(W_{1}\otimes W_{2}, {\cal{D}}_{P(z)}(W))$.
If $V$ is regular, i.e., any weak $V$-module
is a direct sum of irreducible (ordinary) modules, then
we obtain a decomposition of ${\cal{D}}_{P(z)}(W)$ into irreducible 
$V\otimes V$-modules with the corresponding fusion rules as
the multiplicities. In the case that $W=V$, we obtain a decomposition
of ${\cal{D}}_{P(z)}(V)$ of Peter-Weyl type.

This paper is intimately related to [H3, HL0-3] and
we use certain basic principles 
and techniques that have been explored  herein.
In the tensor product theory
for a vertex operator algebra [H3, HL0-3], for {\em two} given $V$-modules 
$W_{1}$ and $W_{2}$,
a $V$-module ``$W_{1}\otimes_{P(z)} W_{2}$'' was constructed such that
for any $V$-module $W$,
a $P(z)$-intertwining map of type ${W'\choose W_{1}W_{2}}$
canonically gives rise to a $V$-homomorphism
{}from ``$W_{1}\otimes _{P(z)}W_{2}$'' to $W'$. In the present work,
by constructing a $V\otimes V$-module ${\cal{D}}_{P(z)}(W)$ for {\em one}
 given $V$-module $W$,
a $P(z)$-intertwining map of type ${W'\choose W_{1}W_{2}}$
becomes a homomorphism in the category of (weak) {\em $V\otimes V$-modules}
instead of the category of {\em $V$-modules}.

In a recent work [FM], Frenkel and Malikov obtained some interesting
results by using Kazhdan-Lusztig tensoring and Harish-Chandra categories.
The bimodule ${\cal{D}}_{P(z)}(W)$ is analogous to the
Harish-Chandra bimodule associated to a module for a Lie algebra.
We hope to study ${\cal{D}}_{P(z)}(V)$ in this aspect
in a future publication.
During this research we have noticed that regular representations
for affine Lie algebras have been studied in [FP] and that
certain results of Peter-Weyl type on Kac-Moody groups have been obtained
in [KP]. Presumably, [FP], [KP] (in the affine case) and the present 
paper are closely related.
We hope to study the connections in a future publication.

In the next paper [Li6] we shall define and study ${\cal{D}}_{Q(z)}(W)$ 
for a given $V$-module $W$ and we shall study the connection between 
the regular representation ${\cal{D}}_{Q(z)}(V)$
and the modular invariance property  of trace functions
and correlation functions ([Z], [DLM3]).

This paper is organized as follows: In Section 2
we review the notion of the contragredient module and 
discuss certain related issues. In Section 3, we 
present $V\otimes V$-module ${\cal{D}}_{P(z)}(W)$. In Section 4,
we present the analogue of the Peter-Weyl theorem. 

{\bf Acknowledgments.} We thank James Lepowsky for many interesting 
discussions on certain subtle issues in formal calculus.
In this research we greatly benefit from studying Huang and 
Lepowsky's papers [H3], [HL0-3] and we adopt many of their 
important viewpoints into this paper.
We would like also to thank IAS for its financial support
and hospitality while most of this work was done during a membership
in the Spring of 1997. 

\newpage
\section{Contragredient modules and some related issues}
This section is preliminary. Here we review some basic notions and facts, 
including the notion of contragredient module and an extension of this.
We shall discuss the notion of right module and prove certain results.

We shall use standard definitions and notions
such as the notions of vertex operator algebra, (irreducible) module,
homomorphism, for which we 
refer the reader to [FHL] and [FLM]. Since we here only deal with the 
representation theory of vertex operator algebras, we fix a vertex 
operator algebra $V$ throughout this paper.
We shall typically use $x,y, x_{1},x_{2}, x_{3},\dots$ for 
mutually commuting formal variables, use letters $u,v$ for elements of 
$V$ and use $w, w_{(1)},w_{(2)},\dots$ for elements of 
$V$-modules. We use ${\Z}$ for the set of all integers, ${\N}$ for the set of
all nonnegative integers and ${\C}$ for the field of complex numbers.
We shall also use certain relatively new notions, which we recall next. 

\bd{dweak} {\em A {\em weak} $V$-module is a vector space $W$ equipped 
with a linear map $Y_{W}$, called the {\em vertex operator map}, from 
$V$ to $({\rm End}W)[[x,x^{-1}]]$ 
such that the {\em truncation condition}:
\begin{eqnarray}\label{etrun}
Y_{W}(v,x)w\in W((x))\;\;\;\mbox{ for }v\in V, \;w\in W;
\end{eqnarray}
the {\em vacuum property}:
\begin{eqnarray}\label{evacuum}
Y_{W}({\bf 1},x)=1\;\;\;(1\mbox{ on the right being the 
identity operator on }W);
\end{eqnarray}
and the {\em Jacobi identity}:
\begin{eqnarray}\label{ejacobi}
& &x_{0}^{-1}\delta\left(\frac{x_{1}-x_{2}}{x_{0}}\right)
Y_{W}(u,x_{1})Y_{W}(v,x_{2})
-x_{0}^{-1}\delta\left(\frac{x_{2}-x_{1}}{-x_{0}}\right)
Y_{W}(v,x_{2})Y_{W}(u,x_{1})\nonumber\\
&=&x_{2}^{-1}\delta\left(\frac{x_{1}-x_{0}}{x_{2}}\right)
Y_{W}(Y(u,x_{0})v,x_{2})
\end{eqnarray}
for $u,v\in V$, hold on $W$.}
\ed

\br{rder}
{\em It was proved in [DLM2] that the $L(-1)$-derivative property:
\begin{eqnarray}
Y_{W}(L(-1)v,x)={d\over dx}Y_{W}(v,x)
\end{eqnarray}
for $v\in V$ holds for a weak $V$-module $W$. Thus a weak $V$-module 
satisfies all the axioms for a $V$-module given 
in [FHL] and [FLM] except for those involving the $L(0)$-grading.}
\er

{\em Throughout this paper, a $V$-module always stands 
for an ordinary $V$-module
unless it is labeled as a weak or generalized  $V$-module.}

Vertex operator algebra $V$ is said to be {\em regular} [DLM2] 
if any weak $V$-module is a direct sum of irreducible (ordinary) 
$V$-modules. If $V$ is regular, it was proved in  [DLM2] that
$V$ has only finitely many inequivalent irreducible modules. 
Examples of regular  vertex operator algebras were given in [DLM2].

It is well known (cf. [B], [FFR], [Li1], [Li4], [MP]) that there is a 
Lie algebra $g(V)$ associated to $V$. More precisely,
\begin{eqnarray}
g(V)=\hat{V}/d\hat{V},
\end{eqnarray}
where 
$$\hat{V}=V\otimes {\C}[t,t^{-1}]\;\;\mbox{and } 
d=L(-1)\otimes 1+1\otimes {d\over dt},$$
with the following Lie bracket:
\begin{eqnarray}\label{ecomm}
[u(m),v(n)]=\sum_{i\ge 0}{m\choose i}(u_{i}v)(m+n-i)
\end{eqnarray}
for $u,v\in V,\; m,n\in {\Z}$, where $u(m)=u\otimes t^{m}+d\hat{V}$.
Furthermore, $g(V)$ is a ${\Z}$-graded Lie algebra where 
$\deg u(m)=\wt u-m-1$ for homogeneous $u\in V$ and for $m\in {\Z}$.

Let $W$ be a $g(V)$-module. A vector $w$ of $W$ is 
said to be {\em restricted} if for every $v\in V$,
$v(m)w=0$ for $m$ sufficiently large. Denote by $W^{res}$ the
space of all the restricted vectors of $W$.
Since for any $u,v\in V$, $u_{i}v\ne 0$ for only finitely many 
nonnegative integer $i$, it follows from (\ref{ecomm}) that 
$W^{res}$ is a submodule of $W$. 
A $g(V)$-module $W$ is said to be {\em restricted} if $W=W^{res}$.
In general, it is clear that $W^{res}$ is
the unique maximal restricted submodule of $W$.
On the other hand, because of the truncation condition
any weak $V$-module $W$ is a restricted 
$g(V)$-module where $v(n)$ is represented 
by $v_{n}$ for $v\in V, \;n\in {\Z}$. 

We now recall the theory of {\em contragredient module}
established in [FHL].
Let $W=\coprod_{h\in {\C}}W_{(h)}$ be a $V$-module and
let $W'=\coprod_{h\in {\C}}W_{(h)}^{*}$, 
the {\em restricted dual} of $W$. For $v\in V,\; w'\in W'$, define
\begin{eqnarray}\label{ey'}
\<Y'(v,x)w',w\>
=\left\<w',Y\left(e^{xL(1)}\left(-x^{-2}\right)^{L(0)}v,x^{-1}\right)w
\right\>
\end{eqnarray}
for $w\in W$.
Then we have the following fundamental result due to Frenkel, 
Huang and Lepowsky ([FHL], Theorem 5.2.1 and Proposition 5.3.1).

\bp{pfhl0}
The pair $(W',Y')$ carries the structure of a $V$-module and 
$(W'',Y'')=(W,Y)$.
\ep

Let $W$ be a weak $V$-module for now. 
For $v\in V$, set (cf. [HL1])
\begin{eqnarray}\label{eyovx}
Y^{o}(v,x)=Y(e^{xL(1)}(-x^{-2})^{L(0)}v,x^{-1})\in (\End\; W)[[x,x^{-1}]].
\end{eqnarray}
(Note that $e^{xL(1)}(-x^{-2})^{L(0)}v\in V[x,x^{-1}]$.)
For instance, 
\begin{eqnarray}
& &Y^{o}({\bf 1},x)=1,\label{yo1}\\
& &Y^{o}(\omega,x)=x^{-4}Y(\omega,x^{-1})\label{eyomega}
\end{eqnarray}
because $L(0){\bf 1}=L(1){\bf 1}=0$, $L(0)\omega=2\omega$ and $L(1)\omega=0$,
where $\omega$ is the Virasoro element of $V$.
Since $e^{xL(1)}(-x^{-2})^{L(0)}v\in V[x,x^{-1}]$ for $v\in V$ and
$Y(u,x^{-1})w\in W((x^{-1}))$ for any $u\in V$, we have
\begin{eqnarray}\label{eyovxw}
Y^{o}(v,x)w\in W((x^{-1}))\;\;\;\mbox{ for }w\in W.
\end{eqnarray}
It was observed in [HL1] that FHL's proof in [FHL] for 
Proposition \ref{pfhl0}
in fact proves the following {\em opposite Jacobi identity}:
\begin{eqnarray}\label{ehl}
& &x_{0}^{-1}\delta\left(\frac{x_{1}-x_{2}}{x_{0}}\right)
Y^{o}(v,x_{2})Y^{o}(u,x_{1})
-x_{0}^{-1}\delta\left(\frac{x_{2}-x_{1}}{-x_{0}}\right)
Y^{o}(u,x_{1})Y^{o}(v,x_{2})\nonumber\\
&=&x_{2}^{-1}\delta\left(\frac{x_{1}-x_{0}}{x_{2}}\right)
Y^{o}(Y(u,x_{0})v,x_{2}).
\end{eqnarray}
As it was mentioned in [HL1] one should think the pair $(W,Y^{o})$ as
a right weak $V$-module. Since we shall discuss 
a little more further, let us make a formal definition here:
A {\em right } weak $V$-module is
a vector space $W$ equipped 
with a linear map $Y_{W}$ from $V$ to $(\End W)[[x,x^{-1}]]$
such that for $u,v\in V$, $w\in W$,
\begin{eqnarray}
& &Y_{W}(v,x)w\in W((x^{-1})),\\
& &Y_{W}({\bf 1},x)=1,
\end{eqnarray}
and such that the opposite Jacobi identity (\ref{ehl}) with $Y^{o}$ being 
replaced by $Y_{W}$ holds.

\br{rrightmodule1}
{\em Similar to Remark \ref{rder}, for a right weak $V$-module $(W,Y_{W})$ we have
\begin{eqnarray}\label{erightL(-1)}
Y_{W}(L(-1)v,x)={d\over dx}Y_{W}(v,x)
\end{eqnarray}
for $v\in V$. From the opposite Jacobi identity we get
\begin{eqnarray}\label{erightbracket}
[Y_{W}(v,x),L(-1)]=Y_{W}(L(-1)v,x).
\end{eqnarray}
Combining (\ref{erightL(-1)}) with (\ref{erightbracket}) we get
\begin{eqnarray}\label{erightbracketder}
[L(-1),Y_{W}(v,x)]=-{d\over dx}Y_{W}(v,x).
\end{eqnarray} }
\er

\br{rrightmodule2}
{\em Given a right weak $V$-module $(W,Y_{W})$, using the opposite Jacobi identity
and Remark \ref{rrightmodule1} we get the Virasoro right module relations on $W$:
\begin{eqnarray}
-[L(m),L(n)]=(m-n)L(m+n)+{1\over 12}(m^{3}-m)\delta_{m+n,0}({\rm rank} V).
\end{eqnarray}}
\er

\br{rrightmodule3}
{\em A {\em right $V$-module} is defined to be a right weak $V$-module
$W$ on which $L(0)$ semisimply acts such that the grading on $W$, given by 
the $L(0)$-eigenspaces, satisfies the same two grading restrictions as 
those in defining the notion of a 
(left) $V$-module. Then $L(-n)$ are locally nilpotent on $W$ for $n\ge 1$.
(Recall that $L(n)$ are locally nilpotent on left modules.)}
\er

Let $U$ be a vector space together with a linear map $Y_{U}$ from
the vertex operator algebra $V$ to $(\End\; U)[[x,x^{-1}]]$,
e.g., $(U,Y_{U})$ is a left or right (weak) $V$-module.
For $v\in V$, we define (cf. (\ref{eyovx}))
\begin{eqnarray}
Y_{U}^{o}(v,x)=Y_{U}(e^{xL(1)}(-x^{-2})^{L(0)}v,x^{-1})
\in (\End\; U)[[x,x^{-1}]].
\end{eqnarray}
Then FHL's proof in [FHL] for 
Proposition \ref{pfhl0} again
in fact proves the following result, which gives
the equivalence between the notion of left $V$-module 
and the notion of right $V$-module:

\bp{pcontragredient}
Let $V$ be a vertex operator algebra and $U$ a vector space 
together with a linear map  $Y_{U}$ from $V$ to $(\End\; U)[[x,x^{-1}]]$. 
Then
$(U,Y_{U})$ is a left (weak) $V$-module if and only if $(U,Y^{o}_{U})$ 
is a right (weak) $V$-module. 
\ep

Notice that 
\begin{eqnarray}
\Res_{x}xY^{o}(\omega,x)=\Res_{x}x^{-3}Y(\omega,x^{-1})=\Res_{x}xY(\omega,x).
\end{eqnarray}
Then $\omega_{1}=L(0)$ is represented by the same operator
on $U$ for $(U,Y_{U})$ and $(U,Y^{o}_{U})$.
Therefore, the assertion in Proposition \ref{pcontragredient} also includes
the grading information.

Let $f(x)=\sum_{n\le N}a_{n}x^{n}\in U((x^{-1}))$ for a vector 
space $U$. Then, for any complex number $z_{0}$,
\begin{eqnarray}
f(x+z_{0}):=\sum_{n\le N}a_{n}(x+z_{0})^{n}
=\sum_{n\le N}\sum_{j\ge 0}{n\choose j}a_{n}z_{0}^{j}x^{n-j}
\end{eqnarray}
exists in $U((x^{-1}))$. 
Let $(W,Y_{W})$ be a right weak 
$V$-module and $z_{0}$ be a complex number.
For $v\in V$, we define $Y_{W}^{(z_{0})}(v,x)\in (\End\;W)[[x,x^{-1}]]$ by
\begin{eqnarray}
Y_{W}^{(z_{0})}(v,x)w=Y_{W}(v,x+z_{0})w=(Y_{W}(v,y)w)|_{y=x+z_{0}}
=e^{z_{0}{d\over dx}}(Y_{W}(v,x)w)
\end{eqnarray}
for $w\in W$.

\bp{ptranslationz}
For any right weak $V$-module $(W,Y_{W})$, the pair
$(W,Y_{W}^{(z_{0})})$ is also a right weak $V$-module.
\ep

\pf Clearly,
\begin{eqnarray}
Y_{W}^{(z_{0})}({\bf 1},x)=1, \;\; Y_{W}^{(z_{0})}(v,x)w\in W((x^{-1}))
\end{eqnarray}
for $v\in V,\; w\in W$.
Then it remains to prove the opposite Jacobi identity.

Let $u,v\in V$. 
Replacing $x_{1}$ and $x_{2}$ by $x_{1}+z_{0}$ and $x_{2}+z_{0}$
in the opposite Jacobi identity (\ref{ehl}) with $Y^{o}$ being 
replaced by $Y_{W}$, respectively, we get
\begin{eqnarray}\label{ehly}
& &x_{0}^{-1}\delta\left(\frac{x_{1}-x_{2}}{x_{0}}\right)
Y_{W}(v,x_{2}+z_{0})Y_{W}(u,x_{1}+z_{0})\nonumber\\
& &-x_{0}^{-1}\delta\left(\frac{x_{2}-x_{1}}{-x_{0}}\right)
Y_{W}(u,x_{1}+z_{0})Y_{W}(v,x_{2}+z_{0})\nonumber\\
&=&(x_{2}+z_{0})^{-1}\delta\left(\frac{x_{1}+z_{0}-x_{0}}{x_{2}+z_{0}}\right)
Y_{W}(Y(u,x_{0})v,x_{2}+z_{0}),
\end{eqnarray}
noting that for $n\in {\Z}$,
\begin{eqnarray}
(x_{1}+z_{0}-x_{2}-z_{0})^{n}=(x_{1}-x_{2})^{n}.
\end{eqnarray}
Notice that
\begin{eqnarray}
(x_{2}+z_{0})^{-1}\delta\left(\frac{x_{1}+z_{0}-x_{0}}{x_{2}+z_{0}}\right)
=e^{z_{0}({\partial\over \partial x_{1}}+{\partial\over \partial x_{2}})}
x_{2}^{-1}\delta\left(\frac{x_{1}-x_{0}}{x_{2}}\right)
=x_{2}^{-1}\delta\left(\frac{x_{1}-x_{0}}{x_{2}}\right)
\end{eqnarray}
because
\begin{eqnarray}
\left({\partial \over\partial x_{1}}+{\partial \over\partial x_{2}}\right)
x_{2}^{-1}\delta\left(\frac{x_{1}-x_{0}}{x_{2}}\right)=0
\end{eqnarray}
(cf. [Li2]). Then we obtain
\begin{eqnarray}\label{ehlz}
& &x_{0}^{-1}\delta\left(\frac{x_{1}-x_{2}}{x_{0}}\right)
Y_{W}(v,x_{2}+z_{0})Y_{W}(u,x_{1}+z_{0})\nonumber\\
& &-x_{0}^{-1}\delta\left(\frac{x_{2}-x_{1}}{-x_{0}}\right)
Y_{W}(u,x_{1}+z_{0})Y_{W}(v,x_{2}+z_{0})\nonumber\\
&=&x_{2}^{-1}\delta\left(\frac{x_{1}-x_{0}}{x_{2}}\right)
Y_{W}(Y(u,x_{0})v,x_{2}+z_{0}).
\end{eqnarray}
(This is also a generalization of the opposite Jacobi identity for $Y_{W}$
with $z_{0}=0$.)
This proves that $Y_{W}^{(z_{0})}$ satisfies the opposite Jacobi identity.
$\;\;\;\;\Box$


\br{rrightmoduleshift}
{\em Let $(W,Y_{W})$ be a right weak $V$-module on which $L(-1)$ 
is locally nilpotent, 
e.g., $(W,Y_{W})$ is a right $V$-module (recall Remark \ref{rrightmodule3}).
In view of (\ref{erightbracketder}), we have
\begin{eqnarray}\label{econnew}
Y^{(z_{0})}_{W}(v,x)=Y_{W}(v,x+z_{0})=e^{-z_{0}L(-1)}Y_{W}(v,x)e^{z_{0}L(-1)}
\end{eqnarray}
for $v\in V$. Then $(W,Y_{W}^{(z_{0})})$ carries the structure of a right weak 
$V$-module, which is the transported structure through the linear isomorphism 
$e^{-z_{0}L(-1)}$ from $(W,Y_{W})$. 
This gives a different proof of Proposition \ref{ptranslationz}
in this special case and it also proves that 
$e^{-z_{0}L(-1)}$ is a $V$-isomorphism from $(W,Y_{W})$ to $(W,Y_{W}^{(z_{0})})$.}
\er

\br{rchangemodule} {\em 
Let $(W,Y_{W})$ be a left (weak) $V$-module and $z_{0}$ be a
complex number. We first get a right (weak) $V$-module $(W,Y^{o}_{W})$
by Proposition \ref{pcontragredient}.
In view of Proposition \ref{ptranslationz},
 $(W,(Y^{o}_{W})^{(z_{0})})$ is again a right (weak) $V$-module.
Then using Proposition \ref{pcontragredient} again
we get a left (weak) $V$-module
$(W,((Y^{o}_{W})^{(z_{0})})^{o})$. 
The explicit expression of the vertex operator map for this left module
is given by
\begin{eqnarray}
& &(Y^{o}_{W})^{(z_{0})})^{o}(v,x)w\nonumber\\
&=&(Y_{W}^{o})^{(z_{0})}(e^{xL(1)}(-x^{-2})^{L(0)}v,x^{-1})w\nonumber\\
&=&Y_{W}^{o}(e^{xL(1)}(-x^{-2})^{L(0)}v,x^{-1}+z_{0})w\nonumber\\
&=&Y_{W}\left(e^{(x^{-1}+z_{0})L(1)}(-(x^{-1}+z_{0})^{-2})^{L(0)}
e^{xL(1)}(-x^{-2})^{L(0)}v,{x\over 1+z_{0}x}\right)w
\end{eqnarray}
for $v\in V,\; w\in W$. From (5.3.3) of [FHL] we have
\begin{eqnarray}
x_{1}^{L(0)}L(1)x_{1}^{-L(0)}=x_{1}^{-1}L(1),
\end{eqnarray}
which immediately gives
\begin{eqnarray}
x_{1}^{L(0)}e^{xL(1)}x_{1}^{-L(0)}=e^{xx_{1}^{-1}L(1)}.
\end{eqnarray}
Using this formula we obtain
\begin{eqnarray}
((Y^{o}_{W})^{(z_{0})})^{o}(v,x)w
=Y_{W}\left(e^{-z_{0}(1+z_{0}x)L(1)}(1+z_{0}x)^{-2L(0)}v,
\frac{x}{1+z_{0}x}\right)w.
\end{eqnarray}
Furthermore, if $L(1)$ local nilpotently acts on $W$,
by the conjugation formula (5.2.38) of [FHL] we have (cf. (\ref{econnew}))
\begin{eqnarray}\label{ecomplex}
((Y^{o}_{W})^{(z_{0})})^{o}(v,x)=
e^{-z_{0}L(1)}Y(v,x)e^{z_{0}L(1)}.
\end{eqnarray}
Then $e^{z_{0}L(1)}$ is a $V$-isomorphism from
the resulted left module to the original left module.
This can be considered as a special case of 
the general change-variable ([Z], [H2]).}
\er

Notice that for a $V$-module $W$, 
$W'$ is considerably smaller than $W^{*}$ if $W$ is 
infinite-dimensional --- This turns out to be the case almost all 
of the times. Now we consider the whole space $W^{*}$.
Clearly, the definition (\ref{ey'}) still makes sense if one replaces
$w'$ (an element of $W'$) by any element of $W^{*}$.
Let $W$ be a weak $V$-module for now and 
let $v\in V,\;\alpha\in W^{*}$. We define 
\begin{eqnarray}
Y^{*}(v,x)\alpha\in W^{*}[[x,x^{-1}]]
\end{eqnarray}
by
\begin{eqnarray}\label{e*1}
\<Y^{*}(v,x)\alpha, w\>=\<\alpha,Y^{o}(v,x)w\>
\end{eqnarray}
for $w\in W$.

\br{rhln} {\em Notice that the symbol $Y^{*}(v,x)$ was used 
{\em differently} in [HL1] where
the notion $Y^{*}(v,x)$ defined in [HL1] is the notion
$Y^{o}(v,x)$ in this paper. In our notations, 
$Y^{*}$ goes with $W^{*}$ just as $Y'$ goes with $W'$, and the notation
$Y^{o}$ indicates that $Y^{o}$ is
an analogue of the classical opposite multiplication.}
\er

By Proposition 5.3.1 of [FHL] we have
\begin{eqnarray}\label{e*2}
\<\alpha,Y(v,x)w\>
=\left\<Y^{*}\left(e^{xL(1)}\left(-x^{-2}\right)^{L(0)}v,x^{-1}\right)\alpha,w
\right\>.
\end{eqnarray}
Since $Y^{o}({\bf 1},x)=1$ (recall (\ref{yo1})), we have
\begin{eqnarray}\label{ey*1}
Y^{*}({\bf 1},x)=1.
\end{eqnarray}
For $v\in V$, we set
\begin{eqnarray}
Y^{*}(v,x)=\sum_{n\in {\Z}}v^{*}_{n}x^{-n-1}.
\end{eqnarray}
When $v=\omega$ (the Virasoro element), we set
\begin{eqnarray}
Y^{*}(\omega,x)=\sum_{n\in {\Z}}L^{*}(n)x^{-n-2},
\end{eqnarray}
i.e., $L^{*}(n)=\omega^{*}_{n+1}$ for $n\in {\Z}$.
Then (by (\ref{eyomega}))
\begin{eqnarray}\label{eL*n}
\< L^{*}(n)\alpha,w\>=\<\alpha, L(-n)w\>
\end{eqnarray}
for $n\in {\Z},\; \alpha\in W^{*},\; w\in W$.


\br{rstar1} {\em By (\ref{eyovxw}) and (\ref{e*1}) we have
\begin{eqnarray}
\<Y^{*}(v,x)\alpha,w\>=\<\alpha,Y^{o}(v,x)w\>\in {\C}((x^{-1}))
\end{eqnarray}
for $w\in W$, i.e.,
\begin{eqnarray}
Y^{*}(v,x)\alpha\in {\rm Hom}(W,{\C}((x^{-1}))).
\end{eqnarray}
In terms of components, for any $w\in W$  we have
\begin{eqnarray}
\<v^{*}_{n}\alpha,w\>=0\;\;\;\mbox{ for }n\;\;\mbox{sufficiently small.}
\end{eqnarray}}
\er

The opposite Jacobi identity (\ref{ehl}) directly gives 
the following Jacobi identity in terms of matrix-coefficients:
\begin{eqnarray}\label{ejacobiy*m}
& &x_{0}^{-1}\delta\left(\frac{x_{1}-x_{2}}{x_{0}}\right)
\<Y^{*}(u,x_{1})Y^{*}(v,x_{2})\alpha,w\>\nonumber\\
& &-x_{0}^{-1}\delta\left(\frac{x_{2}-x_{1}}{-x_{0}}\right)
\<Y^{*}(v,x_{2})Y^{*}(u,x_{1})\alpha,w\>\nonumber\\
&=&x_{2}^{-1}\delta\left(\frac{x_{1}-x_{0}}{x_{2}}\right)
\<Y^{*}(Y(u,x_{0})v,x_{2})\alpha,w\>
\end{eqnarray}
for $u,v\in V,\; \alpha\in W^{*},\;w\in W$.

\br{rproof} {\em From (\ref{ejacobiy*m}) we have
the following ``Jacobi identity''
\begin{eqnarray}\label{ejacobiy*}
& &x_{0}^{-1}\delta\left(\frac{x_{1}-x_{2}}{x_{0}}\right)
Y^{*}(u,x_{1})Y^{*}(v,x_{2})\alpha
-x_{0}^{-1}\delta\left(\frac{x_{2}-x_{1}}{-x_{0}}\right)
Y^{*}(v,x_{2})Y^{*}(u,x_{1})\alpha\nonumber\\
&=&x_{2}^{-1}\delta\left(\frac{x_{1}-x_{0}}{x_{2}}\right)
Y^{*}(Y(u,x_{0})v,x_{2})\alpha.
\end{eqnarray}
It is important to notice that {\em unlike the usual Jacobi identity,
(\ref{ejacobiy*}) is not algebraic} in the sense that the coefficient of
each monomial $x_{0}^{r}x_{1}^{s}x_{2}^{l}$ in each of the three terms 
is in general an infinite sum in terms of $u^{*}_{m}$ and $v^{*}_{n}$
although it is a well defined element of $W^{*}$. 
The main reason is that $Y^{*}(v,x)\alpha$ in general 
involves infinitely many negative powers of $x$.
Because of the failure of the truncation condition,
technically speaking $(W^{*},Y^{*})$ does not carry the structure of 
a weak $V$-module. On the other hand, in addition to the
Jacobi identity (\ref{ejacobiy*}), $(W^{*},Y^{*})$ 
also satisfies the vacuum property $Y^{*}({\bf 1},x)=1$
and the $L(-1)$-derivative property
$Y^{*}(L(-1)v,x)={d\over dx}Y^{*}(v,x)$ for $v\in V$, which,
as proved in [DLM2] (cf. Remark \ref{rder}), follows from
(\ref{ejacobiy*m}), the vacuum property and the fact 
$L(-1)v=v_{-2}{\bf 1}$.}
\er

As an immediate consequence we have (cf. [Li5]):

\bc{clie}
For any weak $V$-module $W$, $W^{*}$ is a $g(V)$-module
where $v(n)$ acts as $v_{n}^{*}$ for $v\in V, n\in {\Z}$.
In particular, the following commutator formula 
(cf. (\ref{ecomm})) holds:
\begin{eqnarray}
[u_{m}^{*},v_{n}^{*}]=\sum_{i\ge 0}{m\choose i}(u_{i}v)_{m+n-i}^{*}
\end{eqnarray}
for $u,v\in V,\; m,n\in {\Z}$. $\;\;\;\;\Box$
\ec

{\em From now on we shall freely use the commutator 
formula on $W^{*}$ without explicit comments.} 

The following definition naturally arises from Remark \ref{rproof} 
(cf. [Li5]):

\bd{dD}
{\em For a weak $V$-module $W$, we define $D(W)=(W^{*})^{res}$, the 
largest restricted $g(V)$-submodule of $W^{*}$. That is, $D(W)$ consists
of each vector $\alpha$ such that for every $v\in V$,
\begin{eqnarray}
Y^{*}(v,x)\alpha\in W^{*}((x)),
\end{eqnarray}
i.e., $\;v^{*}_{n}\alpha=0$ for $n$ sufficiently large.}
\ed

As it has been noticed in [Li5], from (\ref{ehl}) we immediately have:

\bp{pfhl2}
The pair $(D(W),Y^{*})$ carries the structure of 
a weak $V$-module and $D(W)$ is the unique maximal weak $V$-module 
in $W^{*}$ with $Y^{*}$ being the vertex operator map.$\;\;\;\;\Box$
\ep

\br{rdw} {\em 
Let $v\in V, \;\alpha\in D(W)$. Then there exists $r\in {\Z}$ such that
$$ Y^{*}(v,x)\alpha\in x^{r}W^{*}[[x]].$$
Using (\ref{eyovxw}), we get
\begin{eqnarray}
\<Y^{*}(v,x)\alpha,w\>=\<\alpha,Y^{o}(v,x)w\>
\in x^{r}{\C}[[x]]\cap {\C}((x^{-1}))=x^{r}{\C}[x]
\subset {\C}[x,x^{-1}]
\end{eqnarray}
for all $w\in W$. That is, $\<\alpha,Y^{o}(v,x)w\>$ is a 
Laurent polynomial in $x$.
In other words, $\<\alpha,Y^{o}(v,x)w\>$ is a meromorphic 
function on the sphere $CP^{1}$ 
with only two possible poles at $0$ and $\infty$. Furthermore,
when $v$ is fixed with $w$ being free,
the orders of the possible pole at $0$ are uniformally bounded.
Conversely,
if $\alpha\in W^{*}$ satisfies the above mentioned properties, then
$\alpha\in D(W)$.}
\er

If $W=\coprod_{h\in {\C}}W_{(h)}$ is an ordinary $V$-module, 
then $W^{*}=\prod_{h\in {\C}}W_{(h)}^{*}$ is a formal completion of $W'$
and $Y^{*}$ is the natural extension of $Y'$ (cf. [HL1]). 
Let $W$ be a $V$-module and let $\alpha\in W^{*}$. Then from
(\ref{eL*n}), $\alpha$ is an eigenvector of 
$L^{*}(0)$ with eigenvalue $h$ if and only if
$\alpha\in W_{(h)}^{*}\subset W'$. Therefore, $W'$ is the largest
generalized submodule of $D(W)$.
Let $\cal{A}$ be the class of vertex operator algebras satisfying 
the condition that $L(0)$ acts semisimply on every weak module, 
i.e., every weak module is a generalized module. 
Then we immediately have (cf. [Li5]):

\bp{pdsub}
If $V$ is of class ${\cal{A}}$ and $W$ is a $V$-module,
then $D(W)=W'$.$\;\;\;\;\Box$
\ep

Suppose that $V$ contains a regular vertex operator 
subalgebra (with the same Virasoro element) $V^{0}$.
Since any weak $V$-module is a weak $V^{0}$-module, $L(0)$ (the same 
for both $V$ and $V^{0}$) acts semisimply on any weak $V$-module. Then
$V$ is of class $\cal{A}$, so that Proposition \ref{pdsub} applies to $V$.

\br{rdm} 
{\em  Let $(W,Y_{W})$ be a right weak $V$-module.
Define $D(W)$ to be the subspace of $W^{*}$, consisting of 
each $\alpha$ such that $Y^{*}(v,x)\alpha\in W((x^{-1}))$ for $v\in V$.
Then $(D(W),Y^{*})$ is the unique maximal right weak $V$-module 
in $W^{*}$ with $Y^{*}$ being the vertex operator map.}
\er

\newpage

\section{The weak $V\otimes V$-module ${\cal{D}}_{P(z)}(W)$}
In this section, given a weak $V$-module $W$ and
a nonzero complex number $z$, we define a canonical subspace 
${\cal{D}}_{P(z)}(W)$ of 
$W^{*}$, consisting of what we call $P(z)$-linear functionals on $W$.
This space ${\cal{D}}_{P(z)}(W)$ contains $D(W)$ defined in Section 2 as a 
subspace. We prove that ${\cal{D}}_{P(z)}(W)$
has a natural weak module structure for the tensor product vertex 
operator algebra $V\otimes V$.

{\em Throughout this section, $W$ is a weak $V$-module and
$z$ is a nonzero complex number. }
Now we introduce our first key notion.

\bd{dgood} {\em
A linear functional $\alpha$ on $W$ is called
a {\em $P(z)$-linear functional} if
for any $v\in V, \;w\in W$, the formal series $\<\alpha,Y^{o}(v,x)w\>$, 
an infinite series in general, absolutely converges to a rational 
function of $x$ in 
${\C}[x,x^{-1},(x-z)^{-1}]$ in the domain $|x|>|z|$ such that
when $v$ is fixed with $w$ being free, the orders of the possible poles 
at $0$ and $z$ for the associated rational functions are bounded by a fixed integer.}
\ed

All $P(z)$-linear functionals on $W$ clearly form
a vector space, which we denote by ${\cal{D}}_{P(z)}(W)$.

Notice that a rational function in $x$ with only (two) possible 
poles at $0$ and $z$ can be considered as a meromorphic function 
with only three possible poles at $0, \infty$ and $z$ on the 
sphere $CP^{1}$.
(Recall from Section 2 that $\alpha\in D(W)$ if and only if 
for any $v\in V,\; w\in W$, the formal series $\<\alpha,Y^{o}(v,x)w\>$ 
is a Laurent polynomial, i.e., a meromorphic function on the sphere $CP^{1}$
with only two possible poles at $0$ and $\infty$.)

\br{rP(z)} {\em The designation of the notion of $P(z)$-linear 
functional is 
due to a close connection, which will be given in Section 4,
between the notion of ${\cal{D}}_{P(z)}(W)$
and Huang and Lepowsky's notion of $P(z)$-intertwining map defined 
in [H3] and [HL0-3]. We copy the following information about 
$P(z)$  from [HL1]:
let $K$ be the moduli space of spheres with punctures 
and local coordinates vanishing at these punctures. Then
$P(z)$ is the element of $K$  containing $CP^{1}$ with ordered punctures
$\infty, z, 0$ and standard 
local coordinates ${1\over w}, w-z, w$, vanishing at these punctures, 
respectively.
However, in this paper we shall use only the algebraic aspect of vertex operator 
algebras and one can simply treat ${\cal{D}}_{P(z)}(W)$ as a vector space
depending on $z$ and $W$.}
\er

\br{rcondomain}
{\em Note that if a series $a(x)\in {\C}((x^{-1}))$ absolutely converges
to a rational function in ${\C}[x,x^{-1},(x-z)^{-1}]$ in the domain
$|x|>R$ for some real number $R$, then $a(x)$ must
absolutely converge in the domain $|x|>|z|$ to the same rational function.
Then in Definition \ref{dgood}, it is enough to assume that
$\<\alpha,Y^{o}(v,x)w\>$ absolutely converges
to a rational function in ${\C}[x,x^{-1},(x-z)^{-1}]$ in the domain
$|x|>R$ for some real number $R$.}
\er

In the following we give several equivalent definitions of
$P(z)$-linear functional on $W$.

\bl{lequivconds} 
Let $\alpha\in W^{*}$. Then all the following conditions on 
$\alpha$ are equivalent:

(a) $\alpha\in {\cal{D}}_{P(z)}(W)$.

(b) For $v\in V$, there exist nonnegative integers $l$ and $k$ such that
 for every $w\in W$,
\begin{eqnarray}
x^{l}(x-z)^{k}\<\alpha,Y^{o}(v,x)w\>
\end{eqnarray}
is a (finite) polynomial in $x$.

(c) For $v\in V$, there exist nonnegative integers $l$ and $k$ such that
\begin{eqnarray}
x^{l}(x-z)^{k}Y^{*}(v,x)\alpha\in W^{*}[[x]].
\end{eqnarray}

(d) For $v\in V$, there exists a nonnegative integer $k$ such that
\begin{eqnarray}\label{econ1}
(x-z)^{k}Y^{*}(v,x)\alpha\in W^{*}((x)).
\end{eqnarray}
\el

\pf From the definition, $\alpha\in {\cal{D}}_{P(z)}(W)$
if and only if for $v\in V$, there exist $l, k\in {\N}$ 
such that for every $w\in W$, the formal series
$$x^{l}(x-z)^{k}\<\alpha,Y^{o}(v,x)w\>$$
absolutely converges in the domain $|x|>|z|$ to 
a rational function in $x$ without (finite) poles.
Notice that a rational function without (finite) poles is
simply a polynomial and any series that absolutely converges
to a polynomial in the domain $|x|>|z|$ must be a polynomial itself.
Then it follows that (a) and (b) are equivalent.

Conditions (c) and (d) are obviously equivalent.
To finish the proof we shall show that (b) and (c) are equivalent.
Clearly, (b) implies (c). Since $Y^{o}(v,x)w\in W((x^{-1}))$, from
(c) we get
\begin{eqnarray}
x^{l}(x-z)^{k}\<\alpha,Y^{o}(v,x)w\>\in {\C}[[x]]\cap {\C}((x^{-1}))
={\C}[x].
\end{eqnarray}
This proves that (c) implies (b).$\;\;\;\;\Box$

\br{ryvx-1} {\em 
Since $e^{xL(1)}\left(-x^{-2}\right)^{L(0)}v\in V[x,x^{-1}]$ 
for $v\in V$, in view of (\ref{e*1}), (\ref{e*2}) and 
Lemma \ref{lequivconds} (b),
$\alpha\in {\cal{D}}_{P(z)}(W)$ if and only if for $v\in V$,
there exist $r,s\in {\N}$ such that for every $w\in W$,
\begin{eqnarray}
x^{r}(x-z)^{s}\< \alpha,Y(v,x^{-1})w\>\in {\C}[x]. 
\end{eqnarray}}
\er

The following result gives the closeness of ${\cal{D}}_{P(z)}(W)$
under the actions of the component operators $v^{*}_{n}$ of each 
vertex operator $Y^{*}(v,x)$.

\bp{pbasic1} For $u,v\in V, \;\alpha\in {\cal{D}}_{P(z)}(W)$, there exists 
$k\in {\N}$ such that {\em for all $n\in {\Z}$},
\begin{eqnarray}
(x-z)^{k}Y^{*}(u,x)v^{*}_{n}\alpha\in W^{*}((x)).
\end{eqnarray}
Furthermore, for $v\in V, \;n\in {\Z},\;\alpha\in {\cal{D}}_{P(z)}(W)$,
\begin{eqnarray}
v^{*}_{n}\alpha\in {\cal{D}}_{P(z)}(W).
\end{eqnarray}
That is, the space ${\cal{D}}_{P(z)}(W)$ is a $g(V)$-submodule of $W^{*}$.
\ep

\pf Let $u,v\in V,\;n\in {\Z}$. From the commutator formula 
we have
\begin{eqnarray}\label{e3.7}
Y^{*}(u,x)v^{*}_{n}\alpha
=v^{*}_{n}Y^{*}(u,x)\alpha-\sum_{i\in {\N}}{n\choose i}x^{n-i}
Y^{*}(v_{i}u,x)\alpha.
\end{eqnarray}
Since $v_{i}u=0$ for all but finitely many $i\in {\N}$, there 
exists $k\in {\N}$ (independent of $n$) such that
$$(x-z)^{k}Y^{*}(u,x)\alpha,\;\;
(x-z)^{k}Y^{*}(v_{i}u,x)\alpha\in W^{*}((x))$$
for all $i\in {\N}$. Then using (\ref{e3.7}) we have
$$(x-z)^{k}Y^{*}(u,x)v^{*}_{n}\alpha\in W^{*}((x)).$$
It follows immediately from Lemma \ref{lequivconds} that 
$v_{n}^{*}\alpha\in {\cal{D}}_{P(z)}(W)$.
$\;\;\;\;\Box$

Let ${\C}(x)$ be the algebra of rational functions in $x$. 
Define linear maps $\iota_{x;\infty}$ and $\iota_{x;0}$ 
(cf. [FLM], [FHL]) from ${\C}(x)$
to ${\C}[x,x^{-1}]]$ and ${\C}((x))$, respectively, such that 
for $f(x)\in {\C}(x)$, $\iota_{x;\infty}f(x)$ and $\iota_{x;0}f(x)$ are
the Laurent series expansions at $x=\infty$ and at $x=0$, respectively.
In particular,
\begin{eqnarray}
& &\iota_{x;\infty}(x-z_{0})^{r}=
\sum_{i\ge 0}{r\choose i}(-z_{0})^{i}x^{r-i};\\
& &\iota_{x;0}(x-z_{0})^{r}=
\sum_{i\ge 0}{r\choose i}(-z_{0})^{r-i}x^{i}
\end{eqnarray}
for $r\in {\Z},\;z_{0}\in {\C}^{\times}$.
It is clear that $\iota_{x;\infty}$ and $\iota_{x;0}$ are
${\C}[x,x^{-1}]$-linear and one-to-one, and they are algebra homomorphisms.
{}From Definition \ref{dgood}, 
for $v\in V,\; \alpha\in {\cal{D}}_{P(z)}(W),\; w\in W$,
the formal series $\<\alpha,Y^{o}(v,x)w\>$ lies in the range of
$\iota_{x;\infty}$. Then $\iota_{x;\infty}^{-1}\<\alpha,Y^{o}(v,x)w\>$
is well defined. Furthermore,
$\<\alpha,Y^{o}(v,x)w\>$ absolutely converges 
to $\iota_{x;\infty}^{-1}\<\alpha,Y^{o}(v,x)w\>$
in the domain $|x|>|z|$.
Using the defined $\iota$-maps we define a new
vertex operator map $Y^{R}$.

\bd{dr} {\em Let $v\in V,\; \alpha\in {\cal{D}}_{P(z)}(W)$.
Then we define 
\begin{eqnarray}
Y^{R}(v,x)\alpha\in W^{*}[[x,x^{-1}]]
\end{eqnarray}
by
\begin{eqnarray}
\<Y^{R}(v,x)\alpha,w\>&=&\iota_{x;0}\iota_{x;\infty}^{-1}
\<\alpha,Y^{o}(v,x)w\>\label{e3.13def}\\
&=&\iota_{x;0}\iota_{x;\infty}^{-1}
\<Y^{*}(v,x)\alpha,w\>
\end{eqnarray}
for $w\in W$.}
\ed

{}From the definition we have 
\begin{eqnarray}\label{eiota-1three}
\iota_{x;0}^{-1}\<Y^{R}(v,x)\alpha,w\>=\iota_{x;\infty}^{-1}
\<\alpha,Y^{o}(v,x)w\>
=\iota_{x;\infty}^{-1}
\<Y^{*}(v,x)\alpha,w\>.
\end{eqnarray}

Let $v\in V,\; \alpha\in {\cal{D}}_{P(z)}(W)$ and let $l,k$ be 
as in Lemma \ref{lequivconds} (b) or (c). 
Since $x^{l}\iota_{x;\infty}^{-1}\<\alpha,Y^{o}(v,x)w\>$ is analytic at $x=0$,
\begin{eqnarray}
\<x^{l}Y^{R}(v,x)\alpha,w\>=\iota_{x;0}(x^{l}\iota_{x;\infty}^{-1}
\<\alpha,Y^{o}(v,x)w\>)\in {\C}[[x]]
\end{eqnarray}
for every $w\in W$. Hence
\begin{eqnarray}
x^{l}Y^{R}(v,x)\alpha\in W^{*}[[x]].
\end{eqnarray}
This proves the following result:

\bl{ltruncationyr}
For $v\in V,\; \alpha\in {\cal{D}}_{P(z)}(W)$,
\begin{eqnarray}
Y^{R}(v,x)\alpha\in W^{*}((x)).\;\;\;\;\Box
\end{eqnarray}
\el

Let $k\in {\N}$ be as in Lemma \ref{lequivconds} (b).
(Note that we may take the same $k$ for (b), (c) and (d).)
Noticing that
\begin{eqnarray}
(x-z)^{k}\<\alpha,Y^{o}(v,x)w\>\in {\C}[x,x^{-1}]
\end{eqnarray}
and that $\iota_{x;\infty}$ and $\iota_{x;0}$ are
${\C}[x,x^{-1}]$-linear, from (\ref{e3.13def})
we get
\begin{eqnarray}
(x-z)^{k}\<Y^{R}(v,x)\alpha,w\>
=\iota_{x;0}\iota_{x;\infty}^{-1}(x-z)^{k}\<\alpha,Y^{o}(v,x)w\>
=(x-z)^{k}\<\alpha,Y^{o}(v,x)w\>
\end{eqnarray}
for all $w\in W$. Therefore we have proved the following result:

\bl{cbasic} For $v\in V,\; \alpha\in {\cal{D}}_{P(z)}(W)$,
\begin{eqnarray}\label{e3.11}
(x-z)^{k}Y^{R}(v,x)\alpha=(x-z)^{k}Y^{*}(v,x)\alpha
\end{eqnarray}
whenever $(x-z)^{k}Y^{*}(v,x)\alpha\in W^{*}((x))$ for $k\in {\N}$.
Furthermore
\begin{eqnarray}\label{e3.12}
(x-z)^{k}\<Y^{R}(v,x)\alpha,w\>
=(x-z)^{k}\<\alpha,Y^{o}(v,x)w\>
\end{eqnarray}
for all $w\in W$. $\;\;\;\;\Box$
\el

Set 
\begin{eqnarray}
Y^{R}(v,x)=\sum_{n\in {\Z}}v^{R}_{n}x^{-n-1}.
\end{eqnarray}
For $v=\omega$ (the Virasoro element), following the tradition we set
\begin{eqnarray}
Y^{R}(\omega,x)=\sum_{n\in {\Z}}L^{R}(n)x^{-n-2}.
\end{eqnarray}

For now, $v^{R}_{n}$ are only linear maps from ${\cal{D}}_{P(z)}(W)$
to $W^{*}$. Next, we shall show that $v^{R}_{n}$ are in fact
linear endomorphisms of ${\cal{D}}_{P(z)}(W)$. We shall need the following result:

\bl{lbasic5}
Let $U$ be a vector space and let
\begin{eqnarray}
f(x)=\sum_{n\in {\Z}}f_{n}x^{-n-1},\;\;g(x)=\sum_{n\in {\Z}}g_{n}x^{-n-1}
\in U[[x,x^{-1}]].
\end{eqnarray}
Suppose that either $f(x)\in U((x))$ or $f(x)\in U((x^{-1}))$ and that
there exist $k\in {\N}$ and $z\in {\C}^{\times}$ such that
\begin{eqnarray}\label{ebasic52}
(x-z)^{k}f(x)=(x-z)^{k}g(x).
\end{eqnarray}
Then for $n\in {\Z}$,
\begin{eqnarray}
f_{n}\in \mbox{\rm linear span }\{g_{m}\;|\;m\ge n\}
\end{eqnarray}
if $f(x)\in U((x))$ and
\begin{eqnarray}
f_{n}\in \mbox{\rm linear span }\{g_{m}\;|\;m\le n+k\}
\end{eqnarray}
if $f(x)\in U((x^{-1}))$.
\el

\pf Assume $f(x)\in U((x))$.
For $n\in {\Z}$, let $s$ be a positive integer 
such that $x^{n+s}f(x)\in U[[x]]$. Then $\Res_{x}x^{n+j}f(x)=0$ 
for $j\ge s$ and furthermore,
\begin{eqnarray}
\Res_{x}x^{n+i}(x-z)^{k}f(x)=0
\end{eqnarray}
for $i\ge s$. Using this and (\ref{ebasic52}) we get
\begin{eqnarray}\label{e3.30f}
f_{n}&=&\Res_{x}x^{n}f(x)\nonumber\\
&=&\Res_{x}x^{n}(-z+x)^{-k}(x-z)^{k}f(x)\nonumber\\
&=&\Res_{x}\sum_{i\ge 0}{-k\choose i}(-z)^{-k-i}x^{n+i}(x-z)^{k}f(x)\nonumber\\
&=&\Res_{x}\sum_{i= 0}^{s-1}{-k\choose i}(-z)^{-k-i}x^{n+i}(x-z)^{k}f(x)
\nonumber\\
&=&\Res_{x}\sum_{i= 0}^{s-1}{-k\choose i}(-z)^{-k-i}x^{n+i}(x-z)^{k}g(x)
\nonumber\\
&\in& \mbox{linear span }\{g_{n}, g_{n+1},\cdots, g_{n+k+s-1}\}.
\end{eqnarray}

If $f(x)\in U((x^{-1}))$, let $s'$ be a positive integer such that 
$x^{n-s'}f(x)\in x^{-2}U[[x^{-1}]]$, so that
\begin{eqnarray}
\Res_{x}x^{n-k-i}(x-z)^{k}f(x)=0
\end{eqnarray}
for $i\ge s'$. Then use $(x-z)^{-k}$ instead of 
$(-z+x)^{-k}$ in (\ref{e3.30f}) to complete the proof. 
$\;\;\;\;\;\Box$

In view of Lemmas \ref{ltruncationyr} and \ref{lbasic5} we immediately have:

\bc{lbasic2} Let $v\in V,\; \alpha\in {\cal{D}}_{P(z)}(W)$ and let 
$k\in {\N}$ be such that (\ref{e3.11}) holds. Then for $n\in {\Z}$,
\begin{eqnarray}\label{relate}
v^{R}_{n}\alpha
\in \mbox{\rm linear span }\{v^{*}_{m}\alpha\;|\; m\ge n\}.
\end{eqnarray}
\ec

As an immediate consequence of Corollary \ref{lbasic2}, 
Proposition \ref{pbasic1}
and Lemma \ref{ltruncationyr} we have:

\bp{prd} Let $v\in V,\; \alpha\in {\cal{D}}_{P(z)}(W)$. Then
\begin{eqnarray}
Y^{R}(v,x)\alpha\in \left({\cal{D}}_{P(z)}(W)\right)((x)).
\end{eqnarray}
\ep

\br{ranother} {\em Here is a slightly different proof 
for Proposition \ref{prd}:
Let $k\in {\N}$ be such that
$(x-z)^{k}Y^{*}(v,x)\alpha\in W^{*}((x))$. Then
using Proposition \ref{pbasic1}, we get
\begin{eqnarray}
(z-x)^{k}Y^{*}(v,x)\alpha\in W^{*}((x))\cap {\cal{D}}_{P(z)}(W)[[x,x^{-1}]]
={\cal{D}}_{P(z)}(W)((x)),
\end{eqnarray}
so that
\begin{eqnarray}\label{ebrac}
(z-x)^{-k}\left[(z-x)^{k}Y^{*}(v,x)\alpha\right]\;\;\;\mbox{exists}
\end{eqnarray}
in ${\cal{D}}_{P(z)}(W)((x))$. 
Note that in (\ref{ebrac}), {\em we must use 
the left and right brackets} because
$(z-x)^{-k}Y^{*}(v,x)\alpha$ in general does not exist 
in $W^{*}[[x,x^{-1}]]$.
On the other hand, since $Y^{R}(v,x)\alpha\in W^{*}((x))$ 
(Lemma \ref{ltruncationyr}), we have
\begin{eqnarray}
(z-x)^{-k}\left((z-x)^{k}Y^{R}(v,x)\alpha\right)
=\left((z-x)^{-k}(z-x)^{k}\right)Y^{R}(v,x)\alpha=Y^{R}(v,x)\alpha.
\end{eqnarray}
In view of this, multiplying both sides of (\ref{e3.11})
by $(z-x)^{-k}$ we obtain }
\begin{eqnarray}\label{edef2}
Y^{R}(v,x)\alpha=(z-x)^{-k}\left[(z-x)^{k}Y^{*}(v,x)\alpha\right]
\in \left({\cal{D}}_{P(z)}(W)\right)((x)).
\end{eqnarray}
\er

\br{rcontain}
{\em Let $\alpha\in D(W)$ (defined in Section 2). Then by definition 
$$Y^{*}(v,x)\alpha\in W^{*}((x))\;\;\;\mbox{ for every }v\in V,$$
hence $\alpha\in {\cal{D}}_{P(z)}(W)$ 
(with $k=0$ in Lemma \ref{lequivconds} (d)).
This shows that $D(W)\subset {\cal{D}}_{P(z)}(W)$. Furthermore,
by Lemma \ref{cbasic},
$$Y^{R}(v,x)\alpha=Y^{*}(v,x)\alpha\;\;\mbox{ for }v\in V,\;\alpha\in D(W).$$
Therefore, the pair $({\cal{D}}_{P(z)}(W), Y^{R})$ extends the pair 
$(D(W),Y^{*})$.}
\er

By Remark \ref{rcontain}, $D(W)$ is contained in the intersection of all
${\cal{D}}_{P(z)}(W)$ for $z\in {\C}^{\times}$.
Conversely, let 
$$\alpha\in {\cal{D}}_{P(z_{1})}(W)\cap {\cal{D}}_{P(z_{2})}(W),$$
where $z_{1},z_{2}$ are two distinct nonzero complex numbers.
Then, for $v\in V,\; w\in W$, the formal series $\<\alpha,Y^{o}(v,x)w\>$
converges in the domain $|x|>{\rm max}\{|z_{1}|,|z_{2}|\}$ to a rational 
function $h(x)$ such that
$$h(x)\in {\C}[x,x^{-1},(x-z_{1})^{-1}]\cap {\C}[x,x^{-1},(x-z_{2})^{-1}].$$
Consequently, $h(x)\in {\C}[x,x^{-1}]$. By Remark \ref{rdw},
$\alpha\in D(W)$.
Thus we have proved the following simple fact:

\bp{pintersection}
Let $W$ be a weak $V$-module. Then
\begin{eqnarray}
\cap_{z\in {\C}^{\times}}{\cal{D}}_{P(z)}(W)=D(W).\;\;\;\;\Box
\end{eqnarray}
\ep

\br{rbim} {\em Here we consider the existence of products of 
certain series. Let 
\begin{eqnarray}
A(x),\;\; B(x)\in {\rm Hom}(V,V((x))).
\end{eqnarray}
Then just like one of the three main terms in the Jacobi identity,
\begin{eqnarray}
x_{0}^{-1}\delta\left(\frac{x_{1}-x_{2}}{x_{0}}\right)A(x_{1})B(x_{2})
\;\;\;\;\mbox{exists}
\end{eqnarray}
in $({\rm End}\;V)[[x_{0},x_{0}^{-1},x_{1},x_{1}^{-1},x_{2},x_{2}^{-1}]]$. 
Furthermore, for any $m,n\in {\Z}$,
\begin{eqnarray}
x_{0}^{-1}\delta\left(\frac{x_{1}-x_{2}}{x_{0}}\right)
(z+x_{1})^{m}(z+x_{2})^{n}A(x_{1})B(x_{2})
\;\;\;\;\mbox{still exists}
\end{eqnarray}
in $({\rm End}\;V)[[x_{0},x_{0}^{-1}, x_{1},x_{1}^{-1},x_{2},x_{2}^{-1}]]$,
or equivalently, for each $l\in {\Z}$,
\begin{eqnarray}\label{elmnab}
(x_{1}-x_{2})^{l}(z+x_{1})^{m}(z+x_{2})^{n}A(x_{1})B(x_{2})
\;\;\;\;\mbox{exists}
\end{eqnarray}
in $({\rm End}\;V)[[x_{1},x_{1}^{-1},x_{2},x_{2}^{-1}]]$.
Indeed, when the expression in (\ref{elmnab}) is
applied to a vector in $V$  its coefficient of a monomial 
$x_{1}^{r}x_{2}^{s}$ is a finite sum by first
considering the coefficient of $x_{2}^{s}$ and then 
considering the coefficient of $x_{1}^{r}$.
Therefore, one may multiply the expression in (\ref{elmnab})
by $(z+x_{1})^{p}(z+x_{2})^{q}$ for any $p,q\in {\Z}$ and
apply the associativity law ([FLM], Chapter 2).}
\er

As our first key result we have:

\bt{tyR}
The pair $({\cal{D}}_{P(z)}(W),Y^{R})$ carries the structure of
a weak $V$-module.
\et

\pf The truncation condition:
$Y^{R}(v,x)\alpha\in {\cal{D}}_{P(z)}(W)((x))$
for $v\in V,\;\alpha\in {\cal{D}}_{P(z)}(W)$ holds (Proposition \ref{prd})
and $Y^{R}({\bf 1},x)=Y^{*}({\bf 1},x)=1$. 
Then it remains to prove the Jacobi identity.

Let $u,v\in V,\; \alpha\in {\cal{D}}_{P(z)}(W)$. 
By Proposition \ref{pbasic1}, there exists
$k_{0}\in {\N}$ such that for {\em all} $m\in {\Z}$,
\begin{eqnarray}\label{ey*v*m}
(x_{1}-z)^{k_{0}}Y^{*}(u,x_{1})v^{*}_{m}\alpha\in W^{*}((x_{1})),
\end{eqnarray}
hence (recall Lemma \ref{cbasic})
\begin{eqnarray}\label{e3.27}
(x_{1}-z)^{k_{0}}Y^{R}(u,x_{1})Y^{*}(v,x_{2})\alpha
=(x_{1}-z)^{k_{0}}Y^{*}(u,x_{1})Y^{*}(v,x_{2})\alpha.
\end{eqnarray}
Let $k'\ge k_{0}$ be such that
$$(x_{2}-z)^{k'}Y^{R}(v,x_{2})\alpha=
(x_{2}-z)^{k'}Y^{*}(v,x_{2})\alpha\in W^{*}((x_{2})).$$
Then by (\ref{e3.27}),
\begin{eqnarray}\label{eyruv}
& &(x_{2}-z)^{k'}(x_{1}-z)^{k'}
Y^{R}(u,x_{1})Y^{R}(v,x_{2})\alpha\nonumber\\
&=&(x_{2}-z)^{k'}(x_{1}-z)^{k'}Y^{*}(u,x_{1})Y^{*}(v,x_{2})\alpha.
\end{eqnarray}
Similarly, there exists $k''\in {\N}$ such that
\begin{eqnarray}\label{eyrvu}
& &(x_{1}-z)^{k''}(x_{2}-z)^{k''}Y^{R}(v,x_{2})Y^{R}(u,x_{1})\alpha
\nonumber\\
&=&(x_{1}-z)^{k''}(x_{2}-z)^{k''}Y^{*}(v,x_{2})Y^{*}(u,x_{1})\alpha.
\end{eqnarray}
Set $k=k'+k''$. For any $w\in W$, using (\ref{eyruv}), (\ref{eyrvu}) 
and (\ref{ejacobiy*m}) we obtain
\begin{eqnarray}\label{e3.30}
& &x_{0}^{-1}\delta\left(\frac{x_{1}-x_{2}}{x_{0}}\right)
(x_{1}-z)^{k}(x_{2}-z)^{k}\<Y^{R}(u,x_{1})Y^{R}(v,x_{2})\alpha,w\>\nonumber\\
& &-x_{0}^{-1}\delta\left(\frac{x_{2}-x_{1}}{-x_{0}}\right)
(x_{1}-z)^{k}(x_{2}-z)^{k}\<Y^{R}(v,x_{2})Y^{R}(u,x_{1})\alpha,w\>\nonumber\\
&=&x_{0}^{-1}\delta\left(\frac{x_{1}-x_{2}}{x_{0}}\right)
(x_{1}-z)^{k}(x_{2}-z)^{k}\<Y^{*}(u,x_{1})Y^{*}(v,x_{2})\alpha,w\>
\nonumber\\
& &-x_{0}^{-1}\delta\left(\frac{x_{2}-x_{1}}{-x_{0}}\right)
(x_{1}-z)^{k}(x_{2}-z)^{k}\<Y^{*}(v,x_{2})Y^{*}(u,x_{1})\alpha,w\>
\nonumber\\
&=&x_{2}^{-1}\delta\left(\frac{x_{1}-x_{0}}{x_{2}}\right)
(x_{1}-z)^{k}(x_{2}-z)^{k}\<Y^{*}(Y(u,x_{0})v,x_{2})\alpha,w\>.
\end{eqnarray}
Let $m$ be {\em any fixed} integer. Since $u_{n}v\ne 0$ for 
only finitely many $n\ge m$,
there is a nonnegative integer $k_{1}$ ({\em depending} on $m$) 
such that for {\em all} $n\ge m$,
\begin{eqnarray}\label{eunv}
(x_{2}-z)^{k_{1}}Y^{*}(u_{n}v,x_{2})\alpha\in {\cal{D}}_{P(z)}(W)((x_{2})),
\end{eqnarray}
hence (recall Lemma \ref{cbasic})
\begin{eqnarray}
(x_{2}-z)^{k_{1}}Y^{R}(u_{n}v,x_{2})\alpha=
(x_{2}-z)^{k_{1}}Y^{*}(u_{n}v,x_{2})\alpha.
\end{eqnarray}
Then 
\begin{eqnarray}\label{econn1}
& &{\rm Res}_{x_{0}}x_{0}^{m}
x_{2}^{-1}\delta\left(\frac{x_{1}-x_{0}}{x_{2}}\right)
(x_{2}-z)^{k_{1}}\<Y^{R}(Y(u,x_{0})v,x_{2})\alpha,w\>\nonumber\\
&=&{\rm Res}_{x_{0}}x_{0}^{m}
x_{2}^{-1}\delta\left(\frac{x_{1}-x_{0}}{x_{2}}\right)
(x_{2}-z)^{k_{1}}\<Y^{*}(Y(u,x_{0})v,x_{2})\alpha,w\>.
\end{eqnarray}
Choosing $k_{1}\ge k$ and then combining (\ref{e3.30}) with (\ref{econn1})
we obtain 
\begin{eqnarray}\label{ejac}
& &{\rm Res}_{x_{0}}x_{0}^{m}
x_{0}^{-1}\delta\left(\frac{x_{1}-x_{2}}{x_{0}}\right)
(x_{1}-z)^{k_{1}}(x_{2}-z)^{k_{1}}\<Y^{R}(u,x_{1})
Y^{R}(v,x_{2})\alpha,w\>\nonumber\\
& &-{\rm Res}_{x_{0}}x_{0}^{m}
x_{0}^{-1}\delta\left(\frac{x_{2}-x_{1}}{-x_{0}}\right)
(x_{1}-z)^{k_{1}}(x_{2}-z)^{k_{1}}\<Y^{R}(v,x_{2})Y^{R}(u,x_{1})\alpha,w\>
\nonumber\\
&=&{\rm Res}_{x_{0}}x_{0}^{m}
x_{2}^{-1}\delta\left(\frac{x_{1}-x_{0}}{x_{2}}\right)
(x_{1}-z)^{k_{1}}(x_{2}-z)^{k_{1}}\<Y^{R}(Y(u,x_{0})v,x_{2})\alpha,w\>.
\end{eqnarray}
In view of Remark \ref{rbim},
we may multiply both sides of (\ref{ejac}) by 
$(-z+x_{1})^{-k_{1}}(-z+x_{2})^{-k_{1}}$ to get
\begin{eqnarray}\label{ejac2}
& &{\rm Res}_{x_{0}}x_{0}^{m}
x_{0}^{-1}\delta\left(\frac{x_{1}-x_{2}}{x_{0}}\right)
\<Y^{R}(u,x_{1})Y^{R}(v,x_{2})\alpha,w\>
\nonumber\\
& &-{\rm Res}_{x_{0}}x_{0}^{m}
x_{0}^{-1}\delta\left(\frac{x_{2}-x_{1}}{-x_{0}}\right)
\<Y^{R}(v,x_{2})Y^{R}(u,x_{1})\alpha,w\>\nonumber\\
&=&{\rm Res}_{x_{0}}x_{0}^{m}
x_{2}^{-1}\delta\left(\frac{x_{1}-x_{0}}{x_{2}}\right)
\<Y^{R}(Y(u,x_{0})v,x_{2})\alpha,w\>.
\end{eqnarray}
Since 
there is no term depending on $m$ in (\ref{ejac2}) except for
$x_{0}^{m}$ and $m$ is an arbitrarily fixed integer, 
dropping $\Res_{x_{0}}x_{0}^{m}$ from
(\ref{ejac2}) we obtain the desired Jacobi identity for $Y^{R}$.
This concludes the proof. $\;\;\;\;\Box$

Since $Y^{R}$ depends on $z$, 
it is necessary, especially in certain situations, for us to use 
the notation $Y^{R}_{P(z)}$ for $Y^{R}$, and we shall do so from now on.

\br{rcomments}
{\em Notice that in the construction of the left weak $V$-module
$({\cal{D}}_{P(z)}(W),Y_{P(z)}^{R})$ we essentially  used the
right weak $V$-module $(W,Y^{o})$ (recall Definition \ref{dr}). In view of this, let
us denote the left weak $V$-module
$({\cal{D}}_{P(z)}(W),Y_{P(z)}^{R})$ by $({\cal{D}}_{P(z)}(W,Y^{o}),Y_{P(z)}^{R})$.
Now, let $(U,Y_{U})$ be a right weak $V$-module. Then we can go through 
the whole procedure replacing $Y^{o}$ by $Y_{U}$ to obtain a
left weak $V$-module denoted by $({\cal{D}}_{P(z)}(U,Y_{U}),Y_{P(z)}^{R})$.
Or, we can first consider the left weak $V$-module 
$(U,Y^{o}_{U})$ (Proposition \ref{pcontragredient}) and then apply Theorem \ref{tyR}
to get the left weak $V$-module $({\cal{D}}_{P(z)}(U,Y_{U}),Y_{P(z)}^{R})$,
noting that $(Y_{U}^{o})^{o}=Y_{U}$ ([FHL], [HL1]).
In the definitions of ${\cal{D}}_{P(z)}(U,Y_{U})$ and $Y_{P(z)}^{R}$
and in the corresponding results,
 $\<\alpha,Y_{U}(v,x)w\>$ plays the role that has been played by
$\<\alpha,Y^{o}(v,x)w\>$ so far.}
\er

Recall that given a left weak $V$-module $(W,Y_{W})$ and 
a nonzero complex number $z$,
we have a right weak $V$-module  $(W,(Y^{o}_{W})^{(z)})$ 
(Propositions \ref{pcontragredient} and \ref{ptranslationz}).
Using this and the linear maps $\iota_{x;\infty}$ and $\iota_{x;0}$
defined earlier we now define another vertex operator map $Y_{P(z)}^{L}$;
a linear map from $V$ to $(\End\; {\cal{D}}_{P(z)}(W))[[x,x^{-1}]]$.

\bd{dleftaction}
{\em Let $W$ be a left weak $V$-module and $z$ a nonzero complex number.
For $v\in V,\; \alpha\in {\cal{D}}_{P(z)}(W)$, we define
\begin{eqnarray}
Y_{P(z)}^{L}(v,x)\alpha\in W^{*}[[x,x^{-1}]]
\end{eqnarray}
by
\begin{eqnarray}\label{edl}
\<Y^{L}_{P(z)}(v,x)\alpha,w\>
&=&\iota_{x;0}\iota_{x;\infty}^{-1}\<\alpha,(Y^{o})^{(z)}(v,x)w\>\\
&=&\iota_{x;0}\iota_{x;\infty}^{-1}\<\alpha,Y^{o}(v,x+z)w\>
\end{eqnarray}
for $w\in W$.}
\ed

Since $\<\alpha,Y^{o}(v,y)w\>$ absolutely converges in the domain 
$|y|>|z|$ to a 
rational function in ${\C}[y,y^{-1},(y-z)^{-1}]$, the formal series 
$\<\alpha,Y^{o}(v,x+z)w\>$ absolutely converges in the domain $|x|>|z|$, 
$|x+z|>|z|$ to a rational function in ${\C}[x,x^{-1},(x+z)^{-1}]$, 
where the possible poles at $y=0, z$ 
are transformed to $x=-z,0$. In particular, 
$\<\alpha,Y^{o}(v,x+z)w\>$ absolutely converges in the domain $|x|>2|z|$
(to the same rational function).
Thus, in view of  Remark \ref {rcondomain} and
Definitions \ref{dleftaction} and \ref{dr} we immediately  have:

\bl{lkey} Let $W$ be a left weak $V$-module and $z$ a nonzero 
complex number. Then
\begin{eqnarray}
({\cal{D}}_{P(z)}(W),Y^{L}_{P(z)})
=({\cal{D}}_{P(-z)}(W,(Y^{o}_{W})^{(z)}), Y_{P(-z)}^{R}).\;\;\;\;\Box
\end{eqnarray}
\el

In view of Remark \ref{rcomments}, combining Lemma \ref{lkey} with 
Theorem \ref{tyR} we immediately have:

\bp{pyLnew} Let $W$ be a left weak $V$-module and 
$z$ a nonzero complex number.
Then the pair $({\cal{D}}_{P(z)}(W),Y_{P(z)}^{L})$
carries the structure of a weak $V$-module. $\;\;\;\;\Box$
\ep

Following the tradition, for the Virasoro element $\omega$ we set
\begin{eqnarray}\label{eLln}
Y_{P(z)}^{L}(\omega,x)=\sum_{n\in {\Z}}L_{P(z)}^{L}(n)x^{-n-2}.
\end{eqnarray}

Recall the following delta-function properties ([FLM], [FHL], [HL1]):
\begin{eqnarray}
x_{0}^{-1}\delta\left(\frac{x-z}{x_{0}}\right)
=x^{-1}\delta\left(\frac{x_{0}+z}{x}\right),
\;\;\;\;\;
x^{-1}\delta\left(\frac{z+x_{0}}{x}\right)
=z^{-1}\delta\left(\frac{x-x_{0}}{z}\right)
\end{eqnarray}
and
\begin{eqnarray}\label{esub}
& &x_{0}^{-1}\delta\left(\frac{z-x}{x_{0}}\right)
f(x_{0},x)=x_{0}^{-1}\delta\left(\frac{z-x}{x_{0}}\right)f(z-x,x)
\end{eqnarray}
for $f(x_{0},x)\in \Hom (U_{1},U_{2}[x_{0},x_{0}^{-1},x,x^{-1}])$,
where $U_{1}$ and $U_{2}$ are any vector spaces. (This can be made
more general (cf. [FLM]), but this is enough for our purpose.)

Let $f(x)\in {\C}[x,x^{-1},(x-z)^{-1}]$. Then
we have (cf. [FHL], [HL1])
\begin{eqnarray}\label{eabstractjacobi}
& &x_{0}^{-1}\delta\left(\frac{x-z}{x_{0}}\right)
\iota_{x;\infty}f(x)-x_{0}^{-1}\delta\left(\frac{z-x}{-x_{0}}\right)
\iota_{x;0}f(x)\nonumber\\
&=&z^{-1}\delta\left(\frac{x-x_{0}}{z}\right)(\iota_{x_{0};0}f(x_{0}+z)).
\end{eqnarray}
In particular,
\begin{eqnarray}\label{e3delta}
x_{0}^{-1}\delta\left(\frac{x-z}{x_{0}}\right)
-x_{0}^{-1}\delta\left(\frac{z-x}{-x_{0}}\right)
=z^{-1}\delta\left(\frac{x-x_{0}}{z}\right).
\end{eqnarray}

For $v\in V,\; \alpha\in {\cal{D}}_{P(z)}(W),\;w\in W$, from the definition we have
\begin{eqnarray}
\iota_{x;0}^{-1}\<Y^{L}_{P(z)}(v,x)\alpha,w\>
=\left(\iota_{y;\infty}^{-1}\<\alpha,Y^{o}(v,y)w\>\right)|_{y=x+z}.
\end{eqnarray}
Recall that
\begin{eqnarray}
\iota_{x;0}^{-1}\<Y^{R}_{P(z)}(v,x)\alpha,w\>
=\iota_{x;\infty}^{-1}\<\alpha,Y^{o}(v,x)w\>
=\iota_{x;\infty}^{-1}\<Y^{*}(v,x)\alpha,w\>.
\end{eqnarray}
In view of this and (\ref{eabstractjacobi}), we immediately have:

\bp{prelationnew}
Let $v\in V,\; \alpha\in {\cal{D}}_{P(z)}(W)$. Then
\begin{eqnarray}\label{sjacnew}
& &x_{0}^{-1}\delta\left(\frac{x-z}{x_{0}}\right)
Y^{*}(v,x)\alpha-x_{0}^{-1}\delta\left(\frac{z-x}{-x_{0}}\right)
Y_{P(z)}^{R}(v,x)\alpha\nonumber\\
&=&z^{-1}\delta\left(\frac{x-x_{0}}{z}\right)Y_{P(z)}^{L}(v,x_{0})\alpha.
\end{eqnarray}
\ep

The following lemma (cf. Lemma \ref{cbasic})
immediately  follows from Proposition \ref{prelationnew}.

\bl{lleft1} Let $v\in V,\; \alpha\in {\cal{D}}_{P(z)}(W)$. Then 
\begin{eqnarray}
(z+x_{0})^{l}Y_{P(z)}^{L}(v,x_{0})\alpha&=&
\Res_{x}x_{0}^{-1}\delta\left(\frac{x-z}{x_{0}}\right)x^{l}Y^{*}(v,x)\alpha\\
&=&(z+x_{0})^{l}Y^{*}(v,x_{0}+z)\alpha
\end{eqnarray}
if $x^{l}Y_{P(z)}^{R}(v,x)\alpha\in W^{*}[[x]]$ for 
$l\in {\N}$.$\;\;\;\;\Box$
\el

\bp{plrc} For $u,v\in V$, 
\begin{eqnarray}\label{elrrl}
Y_{P(z)}^{L}(u,x_{1})Y_{P(z)}^{R}(v,x_{2})
=Y_{P(z)}^{R}(v,x_{2})Y_{P(z)}^{L}(u,x_{1}),
\end{eqnarray}
acting on ${\cal{D}}_{P(z)}(W)$.
\ep

\pf Let $\alpha\in {\cal{D}}_{P(z)}(W)$ and let $k\in {\N}$ be such that
\begin{eqnarray}
(y-z)^{k}Y_{P(z)}^{R}(v,y)\alpha=(y-z)^{k}Y^{*}(v,y)\alpha.
\end{eqnarray}
By Proposition \ref{pbasic1} and Lemma \ref{cbasic}, we may choose $k$ so large that
$$(y-z)^{k}Y_{P(z)}^{R}(v,y)Y^{*}(u,x)\alpha=(y-z)^{k}Y^{*}(v,y)Y^{*}(u,x)\alpha$$
also holds. Furthermore,
since both $Y_{P(z)}^{R}$ and $Y^{*}$ satisfy the weak commutativity, 
which follows from the commutator formula, we may choose $k$ such that
the following also hold:
\begin{eqnarray}
& &(x-y)^{k}Y^{*}(u,x)Y^{*}(v,y)=(x-y)^{k}Y^{*}(v,y)Y^{*}(u,x),\\
& &(x-y)^{k}Y_{P(z)}^{R}(u,x)Y_{P(z)}^{R}(v,y)
=(x-y)^{k}Y_{P(z)}^{R}(v,y)Y_{P(z)}^{R}(u,x).
\end{eqnarray}
Then using (\ref{sjacnew}) and all the above identities  we obtain
\begin{eqnarray}\label{elrcom1}
& &z^{-1}\delta\left(\frac{x-x_{0}}{z}\right)
(y-z)^{k}(x-y)^{k}Y_{P(z)}^{L}(u,x_{0})Y_{P(z)}^{R}(v,y)\alpha\nonumber\\
&=&x_{0}^{-1}\delta\left(\frac{x-z}{x_{0}}\right)
(y-z)^{k}(x-y)^{k}Y^{*}(u,x)Y_{P(z)}^{R}(v,y)\alpha\nonumber\\
& &-x_{0}^{-1}\delta\left(\frac{z-x}{-x_{0}}\right)
(y-z)^{k}(x-y)^{k}Y_{P(z)}^{R}(u,x)Y_{P(z)}^{R}(v,y)\alpha
\nonumber\\
&=&x_{0}^{-1}\delta\left(\frac{x-z}{x_{0}}\right)
(y-z)^{k}(x-y)^{k}Y^{*}(u,x)Y^{*}(v,y)\alpha
\nonumber\\
& &-x_{0}^{-1}\delta\left(\frac{z-x}{-x_{0}}\right)
(y-z)^{k}(x-y)^{k}Y_{P(z)}^{R}(u,x)Y_{P(z)}^{R}(v,y)\alpha\nonumber\\
&=&x_{0}^{-1}\delta\left(\frac{x-z}{x_{0}}\right)
(y-z)^{k}(x-y)^{k}Y^{*}(v,y)Y^{*}(u,x)\alpha
\nonumber\\
& &-x_{0}^{-1}\delta\left(\frac{z-x}{-x_{0}}\right)
(y-z)^{k}(x-y)^{k}Y_{P(z)}^{R}(v,y)Y_{P(z)}^{R}(u,x)\alpha
\nonumber\\
&=&x_{0}^{-1}\delta\left(\frac{x-z}{x_{0}}\right)
(y-z)^{k}(x-y)^{k}Y_{P(z)}^{R}(v,y)Y^{*}(u,x)\alpha
\nonumber\\
& &-x_{0}^{-1}\delta\left(\frac{z-x}{-x_{0}}\right)
(y-z)^{k}(x-y)^{k}Y_{P(z)}^{R}(v,y)Y_{P(z)}^{R}(u,x)\alpha
\nonumber\\
&=&z^{-1}\delta\left(\frac{x-x_{0}}{z}\right)
(y-z)^{k}(x-y)^{k}Y_{P(z)}^{R}(v,y)Y_{P(z)}^{L}(u,x_{0})\alpha.
\end{eqnarray}
Taking $\Res_{x}$ of (\ref{elrcom1}) and then using (\ref{esub}) we get
\begin{eqnarray}\label{elrcom2}
& &(y-z)^{k}(z+x_{0}-y)^{k}Y_{P(z)}^{L}(u,x_{0})Y_{P(z)}^{R}(v,y)\alpha\nonumber\\
&=&(y-z)^{k}(z+x_{0}-y)^{k}Y_{P(z)}^{R}(v,y)Y_{P(z)}^{L}(u,x_{0})\alpha.
\end{eqnarray}
In view of Remark \ref{rbim} we may
multiply both sides of (\ref{elrcom2}) by $(-z+y)^{-k}(z+x_{0}-y)^{-k}$, 
to get 
$$Y_{P(z)}^{L}(u,x_{0})Y_{P(z)}^{R}(v,y)\alpha
=Y_{P(z)}^{R}(v,y)Y_{P(z)}^{L}(u,x_{0})\alpha.$$
This proves (\ref{elrrl}), $\;\;\;\;\Box$

Next we shall combine the two weak $V$-module structures
$Y_{P(z)}^{L}$ and $Y_{P(z)}^{R}$ on ${\cal{D}}_{P(z)}(W)$
into a weak $V\otimes V$-module structure. 
Recall ([FHL], Proposition 3.7.1) that $V\otimes V$ has a natural 
vertex operator algebra structure, where the vertex operator map 
$Y$ is defined by
\begin{eqnarray}
Y(u\otimes v,x)=Y(u,x)\otimes Y(v,x)\;\;\;\mbox{for }u,v\in V. 
\end{eqnarray}
It was proved ([FHL], Propositions 4.6.1)
that for any $V$-modules $W_{1}$ and $W_{2}$, $W_{1}\otimes W_{2}$
has a natural $V\otimes V$-module structure.

Let $Y_{P(z)}(\cdot,x)$ be the (unique) linear map from $V\otimes V$ to 
$(\End \;{\cal{D}}_{P(z)}(W))[[x,x^{-1}]]$ such that
\begin{eqnarray}
Y_{P(z)}(u\otimes v,x)=Y_{P(z)}^{L}(u,x)Y_{P(z)}^{R}(v,x)
\;\;\;\mbox{ for }u,v\in V.
\end{eqnarray}
Now we present our second main theorem of the paper.

\bt{t2} 
The pair $({\cal{D}}_{P(z)}(W),Y_{P(z)})$ carries the structure 
of a weak $V\otimes V$-module.
\et

\pf The proof of Proposition 4.6.1 of [FHL] in fact proves 
the following result:
If we have two commuting left weak $V$-module structures $Y_{1}$ 
and $Y_{2}$ on a vector space $M$ in the sense that
\begin{eqnarray}
Y_{1}(u,x_{1})Y_{2}(v,x_{2})=Y_{2}(v,x_{2})Y_{1}(u,x_{1})
\end{eqnarray}
for $u,v\in V$, then $Y_{1}\otimes Y_{2}$ gives rise
to a left weak $V\otimes V$-module structure on $M$.
Then it immediately follows from Theorem  \ref{tyR} 
and Propositions \ref{plrc} and \ref{pyLnew}.$\;\;\;\;\;\Box$

\br{rproductsum} {\em Let $W$ be a weak $V$-module.
Then from Theorem \ref{t2}, 
$\coprod_{z\in {\C}^{\times}}{\cal{D}}_{P(z)}(W)$ is
a weak $V\otimes V$-module. On the other hand, since
each ${\cal{D}}_{P(z)}(W)$
is a subspace of $W^{*}$, we consider the sum
\begin{eqnarray}
{\cal{S}}(W)=\sum_{z\in {\C}^{\times}}{\cal{D}}_{P(z)}(W).
\end{eqnarray}
Recall (Proposition \ref{pintersection}) that for distinct 
nonzero complex numbers $z_{1}$ and $z_{2}$, we have
$${\cal{D}}_{P(z_{1})}(W)\cap {\cal{D}}_{P(z_{2})}(W)=D(W).$$
Thus $\sum_{z\in \in {\C}^{\times}}{\cal{D}}_{P(z)}(W)$ 
is not a direct sum. Because of this, we need to consider
the existences of an extension of all $Y^{L}_{Q(z)}$ and an extension of 
all $Y^{R}_{Q(z)}$. It turns out that
all $Y^{R}_{P(z)}$ can be put together 
to give a well defined vertex operator action $Y^{R}$ of $V$ 
on ${\cal{S}}(W)$, where
\begin{eqnarray}
\<Y^{R}(v,x)\alpha,w\>=\iota_{x;0}\iota_{x;\infty}^{-1}\<\alpha,Y^{o}(v,x)w\>
\end{eqnarray}
for $v\in V,\; \alpha\in {\cal{S}}(W),\; w\in W$.
Then
$({\cal{S}}(W),Y^{R})$ carries the structure of a weak $V$-module.
However, all $Y^{L}_{P(z)}$ do not give 
a well defined map from $V$ to ${\cal{S}}(W)[[x,x^{-1}]]$. This is
because for $v\in V,\; \alpha\in D(W), w\in W$, by definition
\begin{eqnarray*}
& &\<Y_{P(z_{1})}^{L}(v,x)\alpha,w\>=\iota_{x;0}\iota_{x;\infty}^{-1}
\<\alpha,Y^{o}(v,x+z_{1})w\>,\\
& &\<Y_{P(z_{2})}^{L}(v,x)\alpha,w\>=\iota_{x;0}\iota_{x;\infty}^{-1}
\<\alpha,Y^{o}(v,x+z_{2})w\>,
\end{eqnarray*}
so that $Y_{P(z_{1})}^{L}(v,x)\alpha$ and $Y_{P(z_{2})}^{L}(v,x)\alpha$
are generally different. Therefore, ${\cal{S}}(W)$ is a weak $V$-module,
but not a weak $V\otimes V$-module. }
\er

\br{rsumrightmodule} {\em Motivated by ${\cal{S}}(W)$ we define
a canonical pair $({\cal{D}}(W),Y^{R})$ as follows:

(a) $\alpha\in {\cal{D}}(W)$ if and only if for every $v\in V$, 
there exists a polynomial $f_{v}(x)$ such that 
\begin{eqnarray}
f_{v}(x)\<\alpha,Y^{o}(v,x)w\>\in {\C}[x]
\end{eqnarray}
for all $w\in W$.

(b) For $v\in V,\; \alpha \in {\cal{D}}(W)$, $Y^{R}(v,x)\alpha$ is defined by
\begin{eqnarray}
\<Y^{R}(v,x)\alpha,w\>=\iota_{x;0}\iota_{x;\infty}^{-1}\<\alpha,Y^{o}(v,x)w\>
\end{eqnarray}
for $w\in W$.

By slightly modifying the arguments for
$({\cal{D}}_{P(z)}(W),Y^{R})$ carrying the structure of a weak 
$V$-module, we see that the pair 
$({\cal{D}}(W),Y^{R})$ carries the structure of a weak $V$-module.
It is clear that ${\cal{S}}(W)$ is a submodule,
but we do not know whether ${\cal{D}}(W)$ coincides with ${\cal{S}}(W)$.}
\er

\section{Decomposition of ${\cal{D}}_{P(z)}(W)$ into irreducible 
$V\otimes V$-modules}

In this section, for $V$-modules $W_{1},W_{2}$ and $W$,
we shall identify a $P(z)$-intertwining map of 
type ${W'\choose W_{1}W_{2}}$ in the sense of [H3] and [HL0-3]
with a $V\otimes V$-homomorphism from $W_{1}\otimes W_{2}$ into
${\cal{D}}_{P(z)}(W)$. By using Huang and Lepowsky's one-to-one 
linear correspondences [HL3]
between the space of intertwining operators and the space of 
$P(z)$-intertwining maps of the same type we
obtain canonical linear isomorphisms from the space of intertwining 
operators of type ${W'\choose W_{1}W_{2}}$ to the space of 
$V\otimes V$-homomorphisms from $W_{1}\otimes W_{2}$ to
${\cal{D}}_{P(z)}(W)$. When $V$ is regular, we obtain a decomposition 
of ${\cal{D}}_{P(z)}(W)$ into irreducible $V\otimes V$-modules.
In the case that $W=V$, we obtain an analogue of Peter-Weyl theorem.

A {\em generalized} $V$-module [HL1] is a weak $V$ module on which $L(0)$ 
acts semisimply. That is, a generalized $V$-module satisfies all the axioms
in the notion of a module except the two grading restrictions on homogeneous spaces.
A {\em lower truncated generalized} $V$-module [H2] is 
a generalized $V$-module that also satisfies the lower truncation condition,
one of the two grading restrictions on homogeneous spaces. That is, the only difference 
between a module and a generalized module is that the homogeneous subspaces of a generalized
module could be infinite-dimensional.

Following [FHL] and [HL1], for a vector space $U$, we set
\begin{eqnarray}
U\{x\}=\left\{A(x)=\sum_{h\in {\C}}a_{h}x^{h}\;|\; a_{h}\in U
\;\;\;\mbox{ for }h\in {\C}\right\}.
\end{eqnarray}

{\em Throughout this section, $W, W_{1}, W_{2}$ and $ W_{3}$ are 
generalized $V$-modules.} Now we recall the definition of an intertwining 
operator from [FHL] and [HL1]:

\bd{dinter} {\em {\em An intertwining operator} of type 
${W_{3}\choose W_{1}W_{2}}$ is a linear map $\cal{Y}$ from 
$W_{1}\otimes W_{2}$ to $W_{3}\{x\}$, or equivalently,
\begin{eqnarray}
& &W_{1}\rightarrow (\mbox{Hom}(W_{2},W_{3}))\{x\}\nonumber\\
& &w\mapsto {\cal{Y}}(w,x)=\sum_{n\in {\C}}w_{n}x^{-n-1}
\;\;\;(\mbox{where  }w_{n}\in {\rm Hom}(W_{2},W_{3})) 
\end{eqnarray}
such that ``all the defining properties of a module action that 
makes sense hold.'' That is, for 
$v\in V,\; w_{(1)}\in W_{1},\;w_{(2)}\in W_{2}$, we have the 
lower truncation condition
\begin{eqnarray}
(w_{(1)})_{n}w_{(2)}=0\;\;\;\mbox{ for }n\;\;\mbox{whose real part is 
sufficiently large};
\end{eqnarray}
the $L(-1)$-derivative property
\begin{eqnarray}\label{eintL(-1)}
{\cal{Y}}(L(-1)w_{(1)},x)={d\over dx}{\cal{Y}}(w_{(1)},x);
\end{eqnarray}
and the following Jacobi identity
\begin{eqnarray}\label{eintjacobi}
& &x_{0}^{-1}\delta\left(\frac{x_{1}-x_{2}}{x_{0}}\right)
Y(v,x_{1}){\cal{Y}}(w_{(1)},x_{2})
-x_{0}^{-1}\delta\left(\frac{x_{2}-x_{1}}{-x_{0}}\right)
{\cal{Y}}(w_{(1)},x_{2})Y(v,x_{1})
\nonumber\\
& &=x_{2}^{-1}\delta\left(\frac{x_{1}-x_{0}}{x_{2}}\right)
{\cal{Y}}(Y(v,x_{0})w_{(1)},x_{2}).
\end{eqnarray}
All intertwining operators of this type clearly form a vector space,
denoted by ${\cal{V}}^{W_{3}}_{W_{1}W_{2}}$.}
\ed

The following result can be found in [HL2] (cf. [FHL]):

\bp{phl2} Let $W_{1},W_{2}$ and $W_{3}$ be (ordinary) $V$-modules.
We have the following (canonical) linear isomorphism relations:
\begin{eqnarray}
{\cal{V}}^{W_{3}}_{W_{1}W_{2}}\cong {\cal{V}}^{W_{3}}_{W_{2}W_{1}},\;\;\;\;\;
{\cal{V}}^{W_{3}}_{W_{1}W_{2}}\cong {\cal{V}}^{W_{2}'}_{W_{1}W_{3}'}.
\end{eqnarray}
\ep

Notice that in Sections 2 and 3, only integral powers of $z$ were involved.
In this section, intertwining operators will play
an important role, so that in general complex powers of $z$ 
will be involved.
Let 
${\cal{Y}}(\cdot,x)$ be an intertwining operator of type 
${W'\choose W_{1}W_{2}}$. 
Since for $w_{(1)}\in W_{1},\; w_{(2)}\in W_{2}$,
 ${\cal{Y}}(w_{(1)},x)w_{(2)}$ 
 involves complex powers of $x$,
to consider the evaluation ``${\cal{Y}}(w_{(1)},z)w_{(2)}$'' 
we need to choose a
branch of the log function as in [HL1]. Following [HL1] 
we choose log $z$ so that
\begin{eqnarray}
\log z =\log |z|+i \arg z\;\;\;\mbox{ with }\;\;0\le \arg z< 2\pi,
\end{eqnarray}
and arbitrary values of the log function will be denoted by
\begin{eqnarray}
l_{p}(z)=\log z+2p\pi i
\end{eqnarray}
for $p\in {\Z}$. 

For a $V$-module $W=\coprod_{h\in {\C}}W_{(h)}$,
following [HL1] we define the formal completion 
\begin{eqnarray}
\overline{W}=\prod_{h\in {\C}}W_{(h)}.
\end{eqnarray}
Each $f\in \overline{W}$ can be viewed as a formal sum 
$\sum_{h\in {\C}}f_{h}$, which is a well defined element
of $(W')^{*}$. In view of this we have
\begin{eqnarray}
\overline{W}=(W')^{*},\;\;\;\overline{W'}=W^{*}.
\end{eqnarray}
Furthermore, as noticed in [HL1] the action of a vertex operator 
$Y(v,x)$ on $W$ 
can be naturally extended to $\overline{W}$. 
Then $Y^{*}(v,x)$ acting on 
$W^{*}$ $(=\overline{W'})$ is the natural 
extension of $Y'(v,x)$ on $W'$. Because of this,
we shall also use $Y$ for $Y^{*}$ (as it was done in [HL1]).

Let $U_{1}$ and $U_{2}$ be vector spaces and let 
$$A(x)=\sum_{h\in {\C}}a_{h}x^{h}\in (\Hom (U_{1},U_{2}))\{x\}$$
be such that for every $u_{1}\in U_{1}$, $a_{h}u_{1}=0$ for all 
but finitely many $h$, 
so that $A(x)u_{1}$ is a finite sum. Then we define (cf. [HL1])
\begin{eqnarray}\label{eaxz}
A(e^{l_{p}(z)})=A(x)|_{x=e^{l_{p}(z)}}:
=\sum_{h\in {\C}}e^{hl_{p}(z)}a_{h},
\end{eqnarray}
which is a well defined element of $\Hom (U_{1},U_{2})$.
In the case that $U_{2}={\C}$ and $U_{1}=W$, a $V$-module, 
we have $A(y)\in W'\{y\}$ and 
\begin{eqnarray}
(f(y)A(y))|_{y=e^{l_{p}(z)}}=f(z)A(e^{l_{p}(z)})
\end{eqnarray}
for $f(y)\in {\C}[y,y^{-1}]$.
If ${\cal{Y}}$ is an intertwining operator of type 
${W'\choose W_{1}W_{2}}$, then (cf. [HL1])
\begin{eqnarray}
{\cal{Y}}(w_{(1)},e^{l_{p}(z)})w_{(2)}
={\cal{Y}}(w_{(1)},x)w_{(2)}|_{x=e^{l_{p}(z)}}\in W^{*}
\end{eqnarray}
for $w_{(1)}\in W_{1},\; w_{(2)}\in W_{2}$.

We recall the following notion of $P(z)$-intertwining map from [HL1]:

\bd{rp(z)im} {\rm [HL1]} {\em Let $W_{1}, W_{2}$ and $W$ be generalized $V$-modules.
{\em A $P(z)$-intertwining map}
of type ${W'\choose W_{1}W_{2}}$ is a linear map
$F$ from $W_{1}\otimes W_{2}$ to $\overline{W'} (=W^{*})$ such that
the following identity holds 
for $v\in V,\; w_{(1)}\in W_{1},\; w_{(2)}\in W_{2}$:
\begin{eqnarray}\label{ejacobimap}
& &x_{0}^{-1}\delta\left(\frac{x_{1}-z}{x_{0}}\right)
Y(v,x_{1})F(w_{(1)}\otimes w_{(2)})
-x_{0}^{-1}\delta\left(\frac{z-x_{1}}{-x_{0}}\right)
F(w_{(1)}\otimes Y(v,x_{1})w_{(2)})\nonumber\\
& &=z^{-1}\delta\left(\frac{x_{1}-x_{0}}{z}\right)
F(Y(v,x_{0})w_{(1)}\otimes w_{(2)}).
\end{eqnarray}
Note that $Y(v,x_{1})$ in the expression 
$Y(v,x_{1})F(w_{(1)}\otimes w_{(2)})$
is really $Y^{*}(v,x_{1})$, the natural extension of $Y'(v,x_{1})$ on $W'$.}
\ed

All $P(z)$-intertwining maps of such type clearly form a vector space,
denoted by ${\cal{M}}[P(z)]_{W_{1}W_{2}}^{W'}$.

For $p\in {\Z}$ and for any intertwining 
operator ${\cal{Y}}$ of type ${W'\choose W_{1}W_{2}}$,
we define a linear map (cf. [HL3])
\begin{eqnarray}
F_{{\cal{Y}},p}^{P(z)}: 
& &W_{1}\otimes W_{2}\rightarrow \overline{W'}\;\;(=W^{*})
\nonumber\\
& &\sum_{j=1}^{n}w^{1j}\otimes w^{2j}\mapsto 
\sum_{j=1}^{n}{\cal{Y}}(w_{(1j)},e^{l_{p}(z)})w_{(2j)}
\end{eqnarray}
for $w_{(1j)}\in W_{1},\;w_{(2j)}\in W_{2}$. 
It is clear that $F_{{\cal{Y}},p}^{P(z)}$ is an	
intertwining map of the same type. Furthermore, we have ([HL3], 
Proposition 12.2):

\bp{phl12.2}{\rm [HL3]} Let $W_{1}, W_{2}$ and $W$ be 
lower truncated generalized $V$-modules,
let $z$ be a nonzero complex number and let $p\in {\Z}$. Then
the correspondence ${\cal{Y}}\mapsto F_{{\cal{Y}},p}^{P(z)}$
is a linear isomorphism from the 
space ${\cal{V}}_{W_{1}W_{2}}^{W'}$ of intertwining operators of type
${W'\choose W_{1}W_{2}}$ onto the space ${\cal{M}}[P(z)]_{W_{1}W_{2}}^{W'}$
of $P(z)$-intertwining maps of type ${W'\choose W_{1}W_{2}}$. 
\ep

Now we present our main result in this section.

\bt{tsub}
Let $W_{1}, W_{2}$ and $W$ be generalized $V$-modules and let
$F$ be a linear map from $W_{1}\otimes W_{2}$ to $\overline{W'}$ $(=W^{*})$.
Then $F$ is a $P(z)$-intertwining map of type
${W' \choose W_{1}W_{2}}$ if and only if 
\begin{eqnarray}
F(W_{1}\otimes W_{2})\subset {\cal{D}}_{P(z)}(W)
\end{eqnarray}
and
$F$ is a $V\otimes V$-homomorphism from $W_{1}\otimes W_{2}$ 
into ${\cal{D}}_{P(z)}(W)$ (a subspace of $W^{*}$). Equivalently,
\begin{eqnarray}
{\cal{M}}[P(z)]^{W'}_{W_{1}W_{2}}={\rm Hom}_{V\otimes V}(W_{1}\otimes W_{2},
{\cal{D}}_{P(z)}(W)).
\end{eqnarray}
In particular,
\begin{eqnarray}
\dim {\cal{M}}[P(z)]^{W'}_{W_{1}W_{2}}
=\dim {\rm Hom}_{V\otimes V}(W_{1}\otimes W_{2},
{\cal{D}}_{P(z)}(W)).
\end{eqnarray}
\et

\pf Suppose that $F$ is a $V\otimes V$-homomorphism from 
$W_{1}\otimes W_{2}$ into ${\cal{D}}_{P(z)}(W)$. 
Let $v\in V, \;w_{(1)}\in W_{1},\;w_{(2)}\in W_{2}$.
Since $F(w_{(1)}\otimes w_{(2)})\in {\cal{D}}_{P(z)}(W)$, using
Proposition \ref{prelationnew} and the fact that $F$ is a 
$V\otimes V$-homomorphism we get
\begin{eqnarray}
& &x_{0}^{-1}\delta\left(\frac{x_{1}-z}{x_{0}}\right)
Y(v,x_{1})F(w_{(1)}\otimes w_{(2)})\nonumber\\
&=&x_{0}^{-1}\delta\left(\frac{z-x_{1}}{-x_{0}}\right)
Y_{P(z)}^{R}(v,x_{1})F(w_{(1)}\otimes w_{(2)})+
z^{-1}\delta\left(\frac{x_{1}-x_{0}}{z}\right)
Y_{P(z)}^{L}(v,x_{0})F(w_{(1)}\otimes w_{(2)})\nonumber\\
&=&x_{0}^{-1}\delta\left(\frac{z-x_{1}}{-x_{0}}\right)
F(w_{(1)}\otimes Y(v,x_{1})w_{(2)})+
z^{-1}\delta\left(\frac{x_{1}-x_{0}}{z}\right)
F(Y(v,x_{0})w_{(1)}\otimes w_{(2)}).
\end{eqnarray}
Then $F$ is a $P(z)$-intertwining map.

Conversely, suppose that $F$ is a $P(z)$-intertwining map.
For any $v\in V, \;w_{(1)}\in W_{1}, \;w_{(2)}\in W_{2}$,
let $k\in {\N}$ such that $x_{0}^{k}Y(v,x_{0})w_{(1)}\in W_{1}[[x_{0}]]$.
Then by taking $\Res_{x_{0}}x_{0}^{k}$ from (\ref{ejacobimap}) we obtain
\begin{eqnarray}\label{ercomm}
(x_{1}-z)^{k}Y(v,x_{1})F(w_{(1)}\otimes w_{(2)})
=(x_{1}-z)^{k}F(w_{(1)}\otimes Y(v,x_{1})w_{(2)}).
\end{eqnarray}
Since $Y(v,x_{1})w_{(2)}\in W_{2}((x_{1}))$, the right-hand side of
(\ref{ercomm}) lies in $W^{*}((x_{1}))$, hence
\begin{eqnarray}\label{e4.25}
(x_{1}-z)^{k}Y(v,x_{1})F(w_{(1)}\otimes w_{(2)})\in W^{*}((x_{1})).
\end{eqnarray}
Therefore $F(w_{(1)}\otimes w_{(2)})\in {\cal{D}}_{P(z)}(W)$ 
by Lemma \ref{lequivconds}. Furthermore, by Lemma \ref{cbasic} we have
\begin{eqnarray}
(x_{1}-z)^{k}Y_{P(z)}^{R}(v,x_{1})F(w_{(1)}\otimes w_{(2)})
=(x_{1}-z)^{k}Y(v,x_{1})F(w_{(1)}\otimes w_{(2)}).
\end{eqnarray}
Combining this with (\ref{ercomm}) we get
\begin{eqnarray}\label{etem}
(x_{1}-z)^{k}Y_{P(z)}^{R}(v,x_{1})F(w_{(1)}\otimes w_{(2)})
=(x_{1}-z)^{k}F(w_{(1)}\otimes Y(v,x_{1})w_{(2)}).
\end{eqnarray}
Because 
$$Y_{P(z)}^{R}(v,x_{1})F(w_{(1)}\otimes w_{(2)}),\;\;\;\; 
F(w_{(1)}\otimes Y(v,x_{1})w_{(2)})\in {\cal{D}}_{P(z)}(W)((x_{1})),$$
we can multiply (\ref{etem}) by $(-z+x_{1})^{-k}$ to obtain 
\begin{eqnarray}\label{erhom}
Y_{P(z)}^{R}(v,x_{1})F(w_{(1)}\otimes w_{(2)})
=F(w_{(1)}\otimes Y(v,x_{1})w_{(2)}).
\end{eqnarray}

Let $l\in {\N}$ be such that $x_{1}^{l}Y(v,x_{1})w_{(2)}\in W_{2}[[x_{1}]]$.
Then by taking $\Res_{x_{1}}x_{1}^{l}$ from (\ref{ejacobimap}) we obtain
\begin{eqnarray}\label{e4.24}
\Res_{x_{1}}x_{0}^{-1}\delta\left(\frac{x_{1}-z}{x_{0}}\right)
x_{1}^{l}Y(v, x_{1})F(w_{(1)}\otimes w_{(2)})=
(x_{0}+z)^{l}F(Y(v,x_{0})w_{(1)}\otimes w_{(2)}).
\end{eqnarray}
Similarly, let $l'\in {\N}$ be such that 
$$x_{1}^{l'}Y^{R}_{P(z)}(v,x_{1})F(w_{(1)}\otimes w_{(2)})\in 
({\cal{D}}_{P(z)}(W))[[x_{1}]].$$
By taking $\Res_{x_{1}}x_{1}^{l'}$ from (\ref{sjacnew}) with 
$\alpha=F(w_{(1)}\otimes w_{(2)})$ we obtain
\begin{eqnarray}\label{e4.252}
\Res_{x_{1}}x_{0}^{-1}\delta\left(\frac{x_{1}-z}{x_{0}}\right)
x_{1}^{l'}Y(v, x_{1})F(w_{(1)}\otimes w_{(2)})=
(x_{0}+z)^{l'}Y_{P(z)}^{L}(v,x_{0})F(w_{(1)}\otimes w_{(2)}).
\end{eqnarray}
Combing (\ref{e4.24}) with (\ref{e4.252}) we get
\begin{eqnarray}\label{eteml}
(x_{0}+z)^{l+l'}Y_{P(z)}^{L}(v,x_{0})F(w_{(1)}\otimes w_{(2)})
=(x_{0}+z)^{l+l'}F(Y(v,x_{0})w_{(1)}\otimes w_{(2)}).
\end{eqnarray}
Again, because 
$$Y_{P(z)}^{L}(v,x_{0})F(w_{(1)}\otimes w_{(2)}),\;\;\;\; 
F(Y(v,x_{0})w_{(1)}\otimes w_{(2)})\in ({\cal{D}}_{P(z)}(W))((x_{0})),$$
we can multiply (\ref{eteml}) by $(z+x_{0})^{-l-l'}$ to obtain
\begin{eqnarray}\label{elhom}
Y_{P(z)}^{L}(v,x_{0})F(w_{(1)}\otimes w_{(2)})
=F(Y(v,x_{0})w_{(1)}\otimes w_{(2)}).
\end{eqnarray}
It follows from (\ref{erhom}) and (\ref{elhom}) that 
$F$ is a $V\otimes V$-homomorphism from $W_{1}\otimes W_{2}$
into ${\cal{D}}_{P(z)}(W)$. $\;\;\;\;\Box$

{}From Theorem \ref{tsub}, the image of any $P(z)$-intertwining map of type 
${W'\choose W_{1}W_{2}}$ is contained in ${\cal{D}}_{P(z)}(W)$. 
Let ${\cal{Y}}$ be an intertwining operator of type 
${W'\choose W_{1}W_{2}}$. It follows from Proposition \ref{phl12.2} and
Theorem \ref{tsub} that
\begin{eqnarray}
{\cal{Y}}(w_{(1)},e^{l_{p}(z)})w_{(2)}\in {\cal{D}}_{P(z)}(W)
\end{eqnarray}
for any $w_{(1)}\in W_{1},\;w_{(2)}\in W_{2}, \;p\in {\Z}$.
Furthermore, the linear map $F_{{\cal{Y}},p}^{P(z)}$ 
is a $V\otimes V$-homomorphism.
In view of this, for $p\in {\Z}$  we obtain a linear map
\begin{eqnarray}
F_{p}^{P(z)}:& & {\cal{V}}^{W'}_{W_{1}W_{2}}\rightarrow 
{\rm Hom}_{V\otimes V}(W_{1}\otimes W_{2}, {\cal{D}}_{P(z)}(W))\nonumber\\
& &{\cal{Y}}\mapsto F_{{\cal{Y}},p}^{P(z)}.
\end{eqnarray}
Combining Proposition \ref{phl12.2} with Theorem \ref{tsub} 
we immediately have:

\bc{c3} Let $W_{1},W_{2}, W$ be lower truncated generalized 
$V$-modules and let $p\in {\Z}$. 
Then $F_{p}^{P(z)}$
is a linear isomorphism from ${\cal{V}}^{W'}_{W_{1}W_{2}}$ onto
${\rm Hom}_{V\otimes V}(W_{1}\otimes W_{2},{\cal{D}}_{P(z)}(W))$.
In particular,
\begin{eqnarray}
\dim {\cal{V}}^{W'}_{W_{1}W_{2}}
=\dim {\rm Hom}_{V\otimes V}(W_{1}\otimes W_{2},
{\cal{D}}_{P(z)}(W)).
\end{eqnarray}
\ec

Recall that ${\cal{D}}_{P(z)}(W)$ is a weak 
$V\otimes V$-module in general. We introduce 
the following notion in terms of ordinary $V\otimes V$-submodules of
${\cal{D}}_{P(z)}(W)$:

\bd{drepre} {\em 
A $P(z)$-linear functional $\alpha$ on $W$ is called a 
{\em $P(z)$-representative 
functional} if $\alpha$ generates an ordinary
$V\otimes V$-submodule of ${\cal{D}}_{P(z)}(W)$, on which 
$L_{P(z)}^{L}(0)$ or $L_{P(z)}^{R}(0)$ semisimply acts. Denote by
$R_{P(z)}(W)$ the space of all $P(z)$-representative functionals 
on $W$.}
\ed

Note that on an ordinary $V\otimes V$-submodule of ${\cal{D}}_{P(z)}(W)$,
$L_{P(z)}^{L}(0)$ semisimply acts if and only if $L_{P(z)}^{R}(0)$ 
semisimply acts.
With this definition, any $P(z)$-intertwining map of 
type ${W'\choose W_{1}W_{2}}$ is a $V\otimes V$-homomorphism from 
$W_{1}\otimes W_{2}$ into $R_{P(z)}(W)$.
Because of the two grading restrictions in the definition of 
the notion of a module, $R_{P(z)}(W)$ is a generalized 
$V\otimes V$-module in general.

For the rest of this section, {\em let $S$ be a (non-canonical) 
complete set of representatives of equivalence classes of 
irreducible $V$-modules}. 
{}From [FHL], if $W_{1}$ and $W_{2}$ are irreducible $V$-modules,
$W_{1}\otimes W_{2}$ is an irreducible $V\otimes V$-module and furthermore,
any irreducible $V\otimes V$-module is isomorphic to such a module.
Then $W_{1}\otimes W_{2}$ for $W_{1},W_{2}\in S$ 
form a complete set of representatives of equivalence classes of 
irreducible $V\otimes V$-modules.
Set
\begin{eqnarray}\label{ehvsw}
H_{V,S}^{P(z)}(W)
=\coprod_{W_{1},W_{2}\in S}{\cal{M}}[P(z)]_{W_{1}W_{2}}^{W'}
\otimes (W_{1}\otimes W_{2}),
\end{eqnarray}
which is viewed as a generalized $V\otimes V$-module with the action 
on the second parts $(W_{1}\otimes W_{2})$.
With Theorem \ref{tsub} we obtain a  $V\otimes V$-homomorphism 
$\Psi_{W}$ from $H_{V,S}^{P(z)}(W)$ into $R_{P(z)}(W)$, defined by
\begin{eqnarray}\label{epi}
\Psi_{W} (F\otimes (w_{(1)}\otimes w_{(2)}))=F(w_{(1)}\otimes w_{(2)})
\end{eqnarray}
for $F\in {\cal{M}}[P(z)]_{W_{1}W_{2}}^{W'},\;w_{(1)}\in W_{1},
\; w_{(2)}\in W_{2}$ with $W_{1},W_{2}\in S$.

As we shall need, we discuss certain properties of modules for a 
tensor product vertex operator algebra. For vertex operator algebras 
$V_{1}$ and $V_{2}$, we naturally consider $V_{1}$ and $V_{2}$ as
subalgebras of $V_{1}\otimes V_{2}$ by identifying $V_{1}$ with 
$V_{1}\otimes {\C}{\bf 1}$ and $V_{2}$ with ${\C}{\bf 1}\otimes V_{2}$,
respectively. In this way, any weak $V_{1}\otimes V_{2}$-module is 
automatically a weak $V_{1}$-module and a weak $V_{2}$-module and
the actions of $V_{1}$ and $V_{2}$ commute. (See [FHL], Section 4.7
for more details.) 
We denote by $L^{1}(n)$ and $L^{2}(n)$ for $n\in {\Z}$ 
the corresponding Virasoro operators. Then $L^{1}(n)+L^{2}(n)$ 
for $n\in {\Z}$ are the Virasoro operators for $V_{1}\otimes V_{2}$.

\bp{phom} Let $V_{1}$ and $V_{2}$ be vertex operator algebras,
let $W$ be a weak $V_{1}\otimes V_{2}$-module,
and let $W_{2}$ be a finitely generated weak $V_{2}$-module,
e.g., $W_{2}$ is an irreducible $V_{2}$-module.
Then
${\rm Hom}_{V_{2}}(W_{2},W)$ is a weak $V_{1}$-module with
\begin{eqnarray}
(Y(v_{(1)},x)f)(w_{(2)})=Y(v_{(1)},x)f(w_{(2)})
\end{eqnarray}
for $v_{(1)}\in V_{1},\; f\in \Hom_{V_{2}}(W_{2},W),\;w_{(2)}\in W_{2}$. 
Furthermore, 
if $W_{2}$ is a $V_{1}$-module and $W$ is a $V_{1}\otimes V_{2}$-module
on which $L^{1}(0)$ or $L^{2}(0)$ semisimply acts, then
${\rm Hom}_{V_{2}}(W_{2},W)$ is a $V_{1}$-module.
\ep

\pf Since the actions of $V_{1}$ and $V_{2}$ on $W$ commute, we easily
see that
\begin{eqnarray}
Y(v_{(1)},x)f\in (\Hom_{V_{2}}(W_{2},W))[[x,x^{-1}]]
\end{eqnarray}
for $v_{(1)}\in V_{1},\; f\in \Hom_{V_{2}}(W_{1},W)$.
Clearly, the vacuum property holds and  the Jacobi identity 
will also hold provided that the truncation condition holds.

Let $P$ be a finite-dimensional subspace of $W_{2}$, which 
generates $W_{2}$ as a weak $V_{2}$-module.
 Let $v_{(1)}\in V_{1},\; f\in {\Hom}_{V_{2}}(W_{2},W)$. 
Since $\dim\; P<\infty$ and
$$(Y(v_{(1)},x)f)(p)=Y(v_{(1)},x)f(p)\in W((x))\;\;\;\mbox{ for every }
p\in P,$$
there exists $r\in {\N}$ such that
\begin{eqnarray}
x^{r}(Y(v_{(1)},x)f)(p)=x^{r}Y(v_{(1)},x)f(p)\in W[[x]]\;\;\;\;
\mbox{ for all }p\in P.
\end{eqnarray}
Furthermore, for $v_{(2)}\in V_{2}$ we have
\begin{eqnarray}\label{e4.37}
x_{1}^{r}(Y(v_{(1)},x_{1})f)(Y(v_{(2)},x_{2})p)
&=&Y(v_{(2)},x_{2})x_{1}^{r}Y(v_{(1)},x_{1})f(p)\nonumber\\
&\in& W[[x_{1},x_{2},x_{2}^{-1}]]\;\;\;\;\mbox{ for all }p\in P.
\end{eqnarray}
Since $P$ generates $W_{2}$ as a $V_{2}$-module, 
by repeatedly using (\ref{e4.37}) we get
\begin{eqnarray}
x^{r}(Y(v_{(1)},x)f)(w_{(2)})\in W[[x]]\;\;\;\mbox{ for all }w_{(2)}\in W_{2}.
\end{eqnarray}
Thus $x^{r}Y(v_{(1)},x)f\in \left({\Hom}_{V_{2}}(W_{2},W)\right)[[x]]$.
This proves the truncation condition, so ${\Hom}_{V_{2}}(W_{2},W)$ is a
weak $V_{1}$-module.

Now, assume that $W_{2}$ is a $V_{2}$-module and $W$ is a 
$V_{1}\otimes V_{2}$-module on which $L^{1}(0)$ or $L^{2}(0)$ 
semisimply acts. Since $L^{1}(0)+L^{2}(0)$ semisimply acts on $W$,
both $L^{1}(0)$ and $L^{2}(0)$ semisimply act on $W$.
Then
\begin{eqnarray}
W=\coprod_{\alpha\in {\C}}W[\alpha],
\end{eqnarray}
where
\begin{eqnarray*}
W[\alpha]=\{ w\in W\;|\; L^{1}(0)w=\alpha w\}.
\end{eqnarray*}
Clearly, each $W[\alpha]$ is a weak $V_{2}$-submodule of $W$.
On the other hand, we may take the generating space $P$ of $W_{2}$
to be a graded subspace with 
\begin{eqnarray}
P=P_{(h_{1})}\oplus \cdots \oplus P_{(h_{m})},
\end{eqnarray}
where $h_{1},\dots,h_{m}$ are finitely many distinct  complex numbers.

Let $f\in {\Hom}_{V_{2}}(W_{2},W)$. Since $\dim f(P)<\infty$, 
there are finitely many 
(distinct) complex numbers $\alpha_{1},\dots,\alpha_{n}$ such that
\begin{eqnarray}
f(P)\subset W[\alpha_{1}]\oplus \cdots\oplus W[\alpha_{n}].
\end{eqnarray}
Since $P$ generates $W_{2}$ under the action of $V_{2}$,
\begin{eqnarray}
f(W_{2})\subset W[\alpha_{1}]\oplus \cdots\oplus W[\alpha_{n}].
\end{eqnarray}
Let $f_{i}$ be the projection of $f$ onto $W[\alpha_{i}]$ for $1\le i\le n$.
Clearly, each $f_{i}$ is a $V_{2}$-homomorphism, hence
$f_{i}\in \Hom_{V_{2}}(W_{2},W)$. 
For $1\le i\le n$, we have
\begin{eqnarray}
L(0)f_{i}=\alpha_{i}f_{i}
\end{eqnarray}
because
$$(L(0)f_{i})(w_{(2)})=L^{1}(0)f_{i}(w_{(2)})=\alpha_{i}f_{i}(w_{(2)})$$
for $w_{(2)}\in W_{2}$.
This proves that $L(0)$ semisimply acts on $\Hom_{V_{2}}(W_{2},W)$.
If $f\in (\Hom_{V_{2}}(W_{2},W))_{(\alpha)}$ and $p\in P_{(h)}$, then
$$(L^{1}(0)+L^{2}(0))f(p)=(L(0)f)(p)+f(L(0)p)=(\alpha+h)f(p).$$
That is, $f(p)\in W_{(\alpha+h)}$.
Since $P$ generates $W_{2}$ under the action of $V_{2}$, through
the restriction map we may identify
$(\Hom_{V_{2}}(W_{2},W))_{(\alpha)}$ as a subspace of
$\sum_{j=1}^{m} \Hom(P_{(h_{j})},W_{(\alpha+h_{j})})$.
Then it follows from the finite-dimensionality of $P$ and
the two grading restrictions on $W$ that
the two grading restrictions on
$\Hom_{V_{2}}(W_{2},W)$ also hold.
Therefore, $\Hom_{V_{2}}(W_{2},W)$ is a $V_{1}$-module.$\;\;\;\;\Box$

We also have the following simple fact:

\bl{lhomdirect}
Let $V$ be a vertex operator algebra, $W_{j}$ $(j\in J)$
be weak $V$-modules and $W$ be a finitely generated weak
$V$-module. Then as a vector space,
\begin{eqnarray}
\Hom_{V}\left(W,\coprod_{j\in J}W_{j}\right)
=\coprod_{j\in J}\Hom_{V}(W,W_{j}).
\end{eqnarray}
\el

\pf If $J$ is a finite index set, it follows easily from 
the proof of the corresponding classical result. In general, 
we shall need the assumption on $W$. It is clear that
\begin{eqnarray}\label{econverse}
\coprod_{j\in J}\Hom_{V}(W,W_{j})\subset 
\Hom_{V}\left(W,\coprod_{j\in J}W_{j}\right).
\end{eqnarray}
For the converse, let $P$ be a finite-dimensional generating 
space of $W$. Then for every $f\in \Hom_{V}(W,\coprod_{j\in J}W_{j})$, 
there are finitely many $j_{1},\dots,j_{n}$ such that
\begin{eqnarray}
f(P)\subset W_{j_{1}}\oplus \cdots \oplus W_{j_{n}}.
\end{eqnarray}
Since $P$ generates $W$, we have
\begin{eqnarray}
f(W)\subset W_{j_{1}}\oplus \cdots \oplus W_{j_{n}}.
\end{eqnarray}
Then
$$f\in \Hom_{V}(W,W_{j_{1}})\oplus \cdots \oplus \Hom_{V}(W,W_{j_{n}}).$$
This proves the converse of (\ref{econverse}), hence
completes the proof. $\;\;\;\;\Box$

Because modules for vertex operator algebras have
finite-dimensional homogeneous subspaces and we work on ${\C}$, 
Schur's Lemma for irreducible modules holds (cf. [FHL], Remark 4.7.1).
Then using the proof of the corresponding classical result
(see for example [BD], Chapter II, Proposition 1.14, where 
Lemma \ref{lhomdirect} is used) we have:

\bp{lfact} Let $V_{1}$ and $V_{2}$ be vertex operator algebras and 
$W$ be a weak $V_{1}\otimes V_{2}$-module such that $W$ is a direct sum of
irreducible (ordinary) $V_{2}$-modules.
Let $S_{2}=\{ M_{i}\;|\;i\in I\}$ 
be a complete set of representatives of equivalence classes of 
irreducible $V_{2}$-modules. For each $i\in I$, define a linear map
\begin{eqnarray}
\Psi_{i} : & &{\rm Hom}_{V_{2}}(M_{i},W)\otimes M_{i}\rightarrow W
\nonumber\\
& & (f,w)\mapsto f(w).
\end{eqnarray}
Then the natural linear map
\begin{eqnarray}
\Psi =(\Psi_{i}): 
\coprod_{i\in I}{\rm Hom}_{V_{2}}(M_{i},W)\otimes M_{i}\rightarrow W
\end{eqnarray}
is a $V_{1}\otimes V_{2}$-isomorphism.
\ep





If every $V\otimes V$-module is completely reducible, then
$R_{P(z)}(W)$, being a sum of $V\otimes V$-modules, is completely reducible,
Then it follows from Theorem \ref{tsub}  and
 Proposition \ref{lfact} with $V_{1}={\C}$ and $V_{2}=V\otimes V$ that
for any $V$-module $W$, $\Psi_{W}$, defined previously, 
is a $V\otimes V$-isomorphism onto $R_{P(z)}(W)$. 

The following is a 
version of a result of [DMZ] about the rationality of the 
tensor product of rational vertex operator algebras.

\bl{ltensorrationality}
Let $V_{1}$ and $V_{2}$ be vertex operator algebras such that
every module is completely reducible. 
Then any $V_{1}\otimes V_{2}$-module $W$, 
on which $L^{1}(0)$ or $L^{2}(0)$ semisimply acts, 
is completely reducible.
\el

\pf Since $L^{1}(0)$ and $L^{2}(0)$ commute, each vector of $W$ is 
a sum of common eigenvectors for $L^{1}(0)$ and $L^{2}(0)$.
Let $w\in W$ be a common eigenvector for $L^{1}(0)$ and $L^{2}(0)$ with
\begin{eqnarray}
L^{1}(0)w=h_{1}w,\;\;L^{2}(0)w=h_{2}w.
\end{eqnarray}
Then $w\in W_{(h_{1}+h_{2})}$, the $(h_{1}+h_{2})$-eigenspace of $L^{1}(0)+L^{2}(0)$. 
Then under the action of 
$V\otimes {\C}$ (or ${\C}\otimes V$), $w$ generates 
an ordinary  $V$-module, which is completely reducible from  the assumption.
Thus $W$ is a completely reducible $V$-module under the action of
$V\otimes {\C}$ (or ${\C}\otimes V$). Let $\{W_{i}\;|\; i\in I\}$ be a 
complete set of non-isomorphic irreducible $V$-submodules of $W$
under the action of $V\otimes {\C}$. 
In view of Proposition  \ref{lfact} we have
$$W=\oplus_{i\in I} \Hom _{V}(W_{i},W)\otimes W_{i}.$$
Furthermore, in view of Proposition \ref{phom}
$\Hom_{V}(W_{i},W)$ is a $V$-module under the action of
${\C}\otimes V$, which is completely reducible. 
Proposition 4.7.2  of [FHL] states that the tensor product of 
irreducible modules for factors is an irreducible module 
for the tensor product vertex operator algebra.
Then it follows that $W$ is a
completely reducible $V\otimes V$-module. $\;\;\;\;\Box$

\br{rrationality} {\em There are certain notions such as regularity [DLM2]
and rationality ([Z], [DLM1], [HL1]), which were defined in terms of
the complete reducibility of certain types of weak modules.
In [HL1], $V$ was defined to be rational if
any $V$-module is completely reducible, if there are only finitely many 
non-isomorphic irreducible $V$-modules and if the fusion rule for any
triple of $V$-modules is finite.} 
\er

If $V$ is rational in the sense of [HL1], then
$H_{V,S}^{P(z)}(W))$ is an (ordinary) $V\otimes V$-module.
Furthermore, it follows from Lemma \ref{ltensorrationality} that
$R_{P(z)}(W)$ is completely reducible. Then using Proposition \ref{lfact}
and Theorem \ref{tsub} we immediately have:

\bt{tw}
Let $V$ be a vertex operator algebra such that every $V$-module 
is completely reducible and let $S$ be a (non-canonical) 
complete set of representatives of equivalence classes of 
irreducible $V$-modules. Then for any $V$-module $W$,
$\Psi_{W}$ is a $V\otimes V$-isomorphism from $H_{V,S}^{P(z)}(W)$ 
onto $R_{P(z)}(W)$. In particular, if $V$ is rational in the sense of [HL1],
the above assertion holds and $R_{P(z)}(W)$ is an (ordinary) $V\otimes V$-module.
$\;\;\;\;\Box$
\et

Let us consider a regular vertex operator algebra $V$, i.e., 
any weak $V$-module is a direct sum of irreducible (ordinary) $V$-modules.
It was proved in [DLM2] that $V\otimes V$ is also regular.
(Proposition \ref{lfact} with $V_{1}=V_{2}=V$ together with 
Proposition 4.7.2 of [FHL] gives another slightly different proof.)
Consequently, 
${\cal{D}}_{P(z)}(W)=R_{P(z)}(W)$ is a direct sum of irreducible 
(ordinary) $V\otimes V$-modules. It follows from Theorem \ref{tw}
that $\Psi_{W}$ is 
a $V\otimes V$-isomorphism onto ${\cal{D}}_{P(z)}(W)$.
Furthermore, it was proved ([Li5], Theorem 3.13) that fusion rules 
for triples of irreducible modules are finite. Thus, ${\cal{D}}_{P(z)}(W)$
is an (ordinary) $V\otimes V$-module.
To summarize we have:

\bt{treg}
Let $V$ be a regular vertex operator algebra 
and let $S$ be a (non-canonical) 
complete set of representatives of equivalence classes of 
irreducible $V$-modules. Then for any $V$-module $W$, 
\begin{eqnarray}
{\cal{D}}_{P(z)}(W)=R_{P(z)}(W),
\end{eqnarray}
${\cal{D}}_{P(z)}(W)$ is 
an (ordinary) $V\otimes V$-module, and the linear map
\begin{eqnarray}
\Psi_{W}:  \coprod_{W_{1},W_{2}\in S}
{\cal{M}}[P(z)]^{W'}_{W_{1}W_{2}}\otimes (W_{1}\otimes W_{2})
\rightarrow {\cal{D}}_{P(z)}(W)
\end{eqnarray}
is a $V\otimes V$-isomorphism.$\;\;\;\;\Box$
\et

Next, we consider the special case that $W=V$.
It was known (cf. [Li0], Remark 2.9) that 
${\cal{V}}^{W_{2}'}_{VW_{1}}\cong {\rm Hom}_{V}(W_{1},W_{2}')$.
Then using Schur's lemma, we find that 
for irreducible $V$-modules  $W_{1}$ and $W_{2}$,
$\dim {\cal{V}}^{W_{2}'}_{VW_{1}}=1$ if 
$W_{1}\cong W_{2}'$, and $0$ otherwise.
In view of Proposition \ref{phl2} we have
$\dim {\cal{V}}^{V'}_{W_{1}W_{1}'}=1$. 

Let $(W,Y_{W})$ be a $V$-module. Then 
$Y_{W}$ is an intertwining operator of type ${W\choose VW}$.
As usual, we use $Y$ for $Y_{W}$. From [FHL], we have an
intertwining operator of type ${W\choose WV}$, defined by
\begin{eqnarray}
Y(w,x)v=e^{xL(-1)}Y(v,-x)w
\end{eqnarray}
for $v\in V,\; w\in W$. For convenience, we refer
the canonical intertwining operators $Y$ of types ${W\choose VW}$ and
${W\choose WV}$ as the {\em standard} intertwining operators.
Furthermore, from [HL2] (see also [FHL]), 
we have an intertwining operator ${\cal{Y}}$
of type ${V'\choose WW'}$, defined by
\begin{eqnarray}\label{eww'v'}
\<{\cal{Y}}(w,x)w',v\>
&=&\<w',Y(e^{xL(1)}e^{\pi i L(0)}x^{-2L(0)}w,x^{-1})v\>\nonumber\\
&=&\<w',e^{x^{-1}L(-1)}Y(v,-x^{-1})e^{xL(1)}e^{\pi i L(0)}x^{-2L(0)}w\>
\nonumber\\
&=&\<e^{x^{-1}L(1)}w',Y(v,-x^{-1})e^{xL(1)}e^{\pi i L(0)}x^{-2L(0)}w\>
\end{eqnarray}
for $w\in W,\;w'\in W',\; v\in V$. 
Then for an irreducible $V$-module $W$ we have
\begin{eqnarray}
{\cal{V}}_{WW'}^{V'}={\C}{\cal{Y}}.
\end{eqnarray}
Then in view of Theorem \ref{tw} we immediately have:

\bt{tv} 
Let $V$ be a vertex operator algebra such that 
every $V$-module is completely reducible
and let $S$ be as before. Define a linear map
\begin{eqnarray}
\Phi^{P(z)}_{V}:  \coprod_{W\in S}W\otimes W' \rightarrow R_{P(z)}(V)
\subset {\cal{D}}_{P(z)}(V)
\end{eqnarray}
by
\begin{eqnarray}
\<\Phi^{P(z)}_{V}(w\otimes w'),v\>
=\<e^{z^{-1}L(1)}w', Y(v,-z^{-1})e^{zL(1)}e^{\pi i L(0)}e^{-2L(0)\log z}w\>
\end{eqnarray}
(cf. (\ref{eww'v'})) for $w\in W,\; w'\in W',\; v\in V$. Then 
 $\Phi^{P(z)}_{V}$ is
a $V\otimes V$-isomorphism onto $R_{P(z)}(V)$. 
In particular, the above assertion holds if $V$ is regular, or rational 
in the sense of [HL1].$\;\;\;\;\Box$
\et

\br{rcft} {\em 
In conformal field theory, the physical Hilbert space is usually 
a direct sum of tensor product of irreducible modules
for the left moving algebra and the corresponding irreducible modules
for the right moving algebra. 
It is very interesting to notice that the $V\otimes V$-module 
$R_{P(z)}(V)$
resembles the physical Hilbert space in conformal field theory.}
\er



{}From Theorem \ref{tw}, for any nonzero complex numbers $z_{1}$ 
and $z_{2}$,
$R_{P(z_{1})}(W)$ and $R_{P(z_{2})}(W)$ are isomorphic generalized 
$V\otimes V$-modules if every $V$-module is completely reducible. 
In the following 
we shall give a canonical $V\otimes V$-isomorphism between 
${\cal{D}}_{P(z_{1})}(W)$ and ${\cal{D}}_{P(z_{2})}(W)$ 
without the complete reducibility assumption on $V$-modules.

We recall the following conjugation formula from [FHL] (Lemma 5.2.3):

\bl{lfhl523}
Let $(M,Y_{M})$ be a generalized $V$-module. Then 
\begin{eqnarray}\label{L(0)con}
x_{0}^{L(0)}Y_{M}(v,x)x_{0}^{-L(0)}=Y_{M}(x_{0}^{L(0)}v,xx_{0})
\end{eqnarray}
for $v\in V$, where $x$ and $x_{0}$ are independent commuting
formal variables. Furthermore, when
replacing $x_{0}$ by a nonzero complex number $z$, we have
\begin{eqnarray}\label{eL(0)conz}
e^{L(0)\log z}Y_{M}(v,x)e^{-L(0)\log z}=Y_{M}(z^{L(0)}v,zx).
\end{eqnarray}
\el

Notice that for a nonzero complex number $z$ and for $p\in {\Z}$,
the linear endomorphism $e^{l_{p}(z)L'(0)}$ of $W'$
is naturally extended
to the endomorphism  $e^{l_{p}(z)L^{*}(0)}$ 
of $W^{*} (=\overline{W'})$.
{}From (\ref{eL(0)conz}) we have
\begin{eqnarray}\label{zL(0)con}
e^{l_{p}(z)L^{*}(0)}Y^{*}(v,x)e^{-l_{p}(z)L^{*}(0)}=Y^{*}(z_{0}^{L(0)}v,zx)
\end{eqnarray}
on $W^{*}$ for $v\in V$. Furthermore, we have (cf. Lemma \ref{lfhl523}):

\bp{pisom}
Let $z, z_{1}$ be nonzero complex numbers and let $W$ be a $V$-module. 
Then the linear automorphism $e^{l_{p}(z)L^{*}(0)}$ of $W^{*}$ maps
${\cal{D}}_{P(z_{1})}(W)$ onto ${\cal{D}}_{P(zz_{1})}(W)$ such that
\begin{eqnarray}
& &e^{l_{p}(z)L^{*}(0)}Y^{R}_{P(z_{1})}(v,x)\alpha
=Y^{R}_{P(zz_{1})}(z^{L(0)}v,zx)e^{l_{p}(z)L^{*}(0)}\alpha\label{eisom1}\\
& &e^{l_{p}(z)L^{*}(0)}Y^{L}_{P(z_{1})}(v,x)\alpha
=Y^{L}_{P(zz_{1})}(z^{L(0)}v,zx)e^{l_{p}(z)L^{*}(0)}\alpha\label{eisom2}
\end{eqnarray}
for $v\in V, \;\alpha\in {\cal{D}}_{P(z_{1})}(W)$.
\ep

\pf For $v\in V,\;\alpha\in {\cal{D}}_{P(z_{1})}(W)$, let
$k\in {\N}$ be such that 
\begin{eqnarray}\label{exz1}
(x-z_{1})^{k}Y^{*}(v,x)\alpha\in W^{*}((x)),
\end{eqnarray}
hence (by Lemma \ref{cbasic}),
\begin{eqnarray}\label{e4.40}
(x-z_{1})^{k}Y^{R}_{P(z_{1})}(v,x)\alpha=(x-z_{1})^{k}Y^{*}(v,x)\alpha.
\end{eqnarray}
Without losing much generality we may assume that $v$ is homogeneous.
Using (\ref{zL(0)con}) and (\ref{exz1}) we get
\begin{eqnarray}
& &(x-zz_{1})^{k}Y^{*}(v,x)e^{l_{p}(z)L^{*}(0)}\alpha\nonumber\\
&=&(x-zz_{1})^{k}e^{l_{p}(z)L^{*}(0)}Y^{*}(z^{-L(0)}v,z^{-1}x)\alpha
\nonumber\\
&=&z^{k-\wt v}(z^{-1}x-z_{1})^{k}e^{l_{p}(z)L^{*}(0)}Y^{*}(v,z^{-1}x)\alpha
\nonumber\\
&\in &W^{*}((x)).
\end{eqnarray}
By Lemma \ref{lequivconds}, 
$e^{l_{p}(z)L^{*}(0)}\alpha\in {\cal{D}}_{P(zz_{1})}(W)$.
Furthermore, by Lemma \ref{cbasic},
\begin{eqnarray}\label{e4.42}
(x-zz_{1})^{k}Y^{R}_{P(zz_{1})}(v,x)e^{l_{p}(z)L^{*}(0)}\alpha
=(x-zz_{1})^{k}Y^{*}(v,x)e^{l_{p}(z)L^{*}(0)}\alpha.
\end{eqnarray}
Using (\ref{e4.40}), (\ref{zL(0)con}) and (\ref{e4.42}) we get
\begin{eqnarray}
& &z^{k}(x-z_{1})^{k}e^{l_{p}(z)L^{*}(0)}Y^{R}_{P(z_{1})}(v,x)\alpha
\nonumber\\
&=&z^{k}(x-z_{1})^{k}e^{l_{p}(z)L^{*}(0)}Y^{*}(v,x)\alpha\nonumber\\
&=&(zx-zz_{1})^{k}Y^{*}(z^{L(0)}v,zx)e^{l_{p}(z)L^{*}(0)}\alpha\nonumber\\
&=&(zx-zz_{1})^{k}Y^{R}_{P(zz_{1})}(z^{L(0)}v,zx)e^{l_{p}(z)L^{*}(0)}\alpha.
\end{eqnarray}
Then multiplying by $z^{-k}(-z_{1}+x)^{-k}$ we obtain (\ref{eisom1}).

Since $\alpha\in {\cal{D}}_{P(z_{1})}(W)$ and 
$e^{l_{p}(z)L^{*}(0)}\alpha\in {\cal{D}}_{P(zz_{1})}(W)$,
by Lemma \ref{lleft1}
there exists $l\in {\N}$ such that
\begin{eqnarray}
& &(x_{0}+z_{1})^{l}Y^{L}_{P(z_{1})}(v,x_{0})\alpha
=\Res_{x_{1}}x_{0}^{-1}\delta\left(\frac{x_{1}-z_{1}}{x_{0}}\right)
(x_{0}+z_{1})^{l}Y^{*}(v,x_{1})\alpha,\label{e4.46}\\
& &(x_{0}+zz_{1})^{l}Y^{L}_{P(zz_{1})}(v,x_{0})e^{l_{p}(z)L^{*}(0)}
\alpha\nonumber\\
&=&\Res_{x_{1}}x_{0}^{-1}\delta\left(\frac{x_{1}-zz_{1}}{x_{0}}\right)
(x_{0}+zz_{1})^{l}Y^{*}(v,x_{1})e^{l_{p}(z)L^{*}(0)}\alpha.\label{e4.47}
\end{eqnarray}
Using (\ref{e4.46}), (\ref{zL(0)con}) and (\ref{e4.47}) we get
\begin{eqnarray}
& &z^{l}(x_{0}+z_{1})^{l}e^{l_{p}(z)L^{*}(0)}Y^{L}_{P(z_{1})}(v,x_{0})\alpha
\nonumber\\
&=&\Res_{x_{1}}x_{0}^{-1}\delta\left(\frac{x_{1}-z_{1}}{x_{0}}\right)
z^{l}(x_{0}+z_{1})^{l}e^{l_{p}(z)L^{*}(0)}Y^{*}(v,x_{1})\alpha
\nonumber\\
&=&\Res_{x_{1}}x_{0}^{-1}\delta\left(\frac{x_{1}-z_{1}}{x_{0}}\right)
z^{l}(x_{0}+z_{1})^{l}Y^{*}(z^{L(0)}v,zx_{1})e^{l_{p}(z)L^{*}(0)}\alpha
\nonumber\\
&=&\Res_{x_{1}}z^{-1}x_{0}^{-1}
\delta\left(\frac{z^{-1}x_{1}-z_{1}}{x_{0}}\right)
z^{l}(x_{0}+z_{1})^{l}Y^{*}(z^{L(0)}v,x_{1})e^{l_{p}(z)L^{*}(0)}\alpha
\nonumber\\
&=&\Res_{x_{1}}(zx_{0})^{-1}
\delta\left(\frac{x_{1}-zz_{1}}{zx_{0}}\right)
z^{l}(x_{0}+z_{1})^{l}Y^{*}(z^{L(0)}v,x_{1})e^{l_{p}(z)L^{*}(0)}\alpha
\nonumber\\
&=&z^{l}(x_{0}+z_{1})^{l}
Y^{L}_{P(zz_{1})}(z^{L(0)}v,zx_{0})e^{l_{p}(z)L^{*}(0)}\alpha.
\end{eqnarray}
Multiplying by $z^{-l}(z_{1}+x_{0})^{-l}$ we obtain (\ref{eisom2}).
$\;\;\;\;\Box$

By setting $v=\omega$ (the Virasoro element) in (\ref{eisom1}) 
and (\ref{eisom2}) we obtain
\begin{eqnarray}
e^{l_{p}(z)L^{*}(0)}L^{R}_{P(z_{1})}(0)
&=&L^{R}_{P(zz_{1})}(0)e^{l_{p}(z)L^{*}(0)},\\
e^{l_{p}(z)L^{*}(0)}L^{L}_{P(z_{1})}(0)
&=&L^{L}_{P(zz_{1})}(0)e^{l_{p}(z)L^{*}(0)}.
\end{eqnarray}
Consequently, $e^{l_{p}(z)L^{*}(0)}$ preserves the weight subspaces.
By (\ref{eisom1}) and (\ref{eisom2}) again,
$e^{l_{p}(z)L^{*}(0)}$ maps an ordinary $V\otimes V$-submodule of 
${\cal{D}}_{P(z_{1})}(W)$
to an ordinary $V\otimes V$-submodule of ${\cal{D}}_{P(z_{2})}(W)$. 
Therefore we have:

\bc{crz12} The linear map $e^{l_{p}(z)L^{*}(0)}$ maps $R_{P(z_{1})}(W)$
onto $R_{P(zz_{1})}(W)$. $\;\;\;\;\Box$
\ec

To achieve our goal we also need the following fact 
(cf. Remark \ref{rchangemodule}):

\bl{lfact1}
Let $(M,Y_{M})$ be a generalized $V$-module (for now). Define
$$Y_{M}^{z}(v,x)=Y_{M}(z^{L(0)}v,zx)$$
for $v\in V$. Then $(M,Y_{M}^{z}(\cdot,x))$ is also
a generalized $V$-module and $e^{l_{p}(z)L(0)}$ is a $V$-isomorphism from
$(M,Y_{M}(\cdot,x))$ onto $(M,Y_{M}^{z}(\cdot,x))$.
\el

\pf Using the conjugation formula (\ref{L(0)con}) for $Y_{M}$, we get
$$e^{L(0)\log (z)}Y_{M}(v,x)e^{-L(0)\log (z)}
=Y_{M}(z^{L(0)}v,zx)=Y_{M}^{z}(v,x).$$
Then $(M,Y_{M}^{z}(\cdot,x))$ carries the structure of a 
generalized $V$-module, which is the
transported structure from $(M,Y_{M}(\cdot,x))$ through 
the linear isomorphism
$e^{l_{p}(z)L(0)}$, and furthermore, $e^{l_{p}(z)L(0)}$
is a $V$-isomorphism from $(M,Y_{M})$ onto
$(M,Y_{M}^{z})$.$\;\;\;\;\Box$

Now we are ready to prove our last result:

\bp{piso}
Let $W$ be a $V$-module, $z, z_{1}$ be nonzero complex numbers and 
let $p\in {\Z}$. Then the linear map
\begin{eqnarray}\label{eisopz}
\sigma_{(p,z,z_{1})}:=
e^{-l_{p}(z)L^{L}_{P(zz_{1})}(0)}e^{-l_{p}(z)L^{R}_{P(zz_{1})}(0)}
e^{l_{p}(z)L^{*}(0)}
\end{eqnarray}
is a $V\otimes V$-isomorphism from 
$R_{P(z_{1})}(W)$ onto $R_{P(zz_{1})}(W)$.
\ep

{\bf Proof.} By definition,
\begin{eqnarray}
Y_{P(z')}(u\otimes v,x)=Y^{L}_{P(z')}(u,x)Y^{R}_{P(z')}(v,x)
\end{eqnarray}
for $u,v\in V$ and for any nonzero complex number $z'$.
Using Corollary \ref{crz12}, Proposition \ref{pisom}
and Lemma \ref{lfact1} we get
\begin{eqnarray}
\sigma_{(p,z,z_{1})}Y_{P(z_{1})}(u\otimes v,x)
&=&\sigma_{(p,z,z_{1})}Y_{P(z_{1})}^{L}(u, x)Y_{P(z_{1})}^{R}(v,x)
\nonumber\\
&=&Y_{P(zz_{1})}^{L}(u,x)Y_{P(zz_{1})}^{R}(v,x)\sigma_{(p,z,z_{1})}
\nonumber\\
&=&Y_{P(zz_{1})}(u\otimes v,x)\sigma_{(p,z,z_{1})}.
\;\;\;\;\Box
\end{eqnarray}

\end{document}